\documentclass[12pt,letterpaper]{amsart}
\usepackage{fullpage}
\usepackage{hyperref}
\usepackage{amssymb}
\usepackage{enumerate}
\usepackage{amsmath}
\usepackage{color}
\usepackage{amsthm}
\usepackage{amssymb}
\usepackage{mathtools}
\usepackage{color}
\usepackage{float}
\usepackage{tikz}
\usepackage{tikz-cd}
\usepackage{wrapfig}
\usetikzlibrary{shapes,decorations}
\usepackage[utf8]{inputenc}
\usepackage{enumitem}
\usepackage[english]{babel}
\usepackage{amsmath,calligra,mathrsfs}
\usepackage{csquotes}
\usepackage[capitalize,noabbrev]{cleveref}
\usepackage[backend=bibtex,maxbibnames=99,sorting=anyt]{biblatex}
\renewbibmacro{in:}{}
\bibliography{References}

\newtheorem{theorem}{Theorem}[section]
\newtheorem*{theorem*}{Theorem}

\newtheorem{proposition}[theorem]{Proposition}
\newtheorem*{proposition*}{Proposition}

\newtheorem{corollary}[theorem]{Corollary}

\theoremstyle{definition}
\newtheorem{definition}[theorem]{Definition}

\theoremstyle{remark}
\newtheorem{remark}[theorem]{Remark}

\newtheorem{lemma}[theorem]{Lemma}

\newtheorem{example}[theorem]{Example}

\newcommand{\cat}{\mathbf}
\newcommand{\floor}[1]{\left\lfloor #1 \right\rfloor}

\newcommand{\id}[1]{\textnormal{id}_{#1}}

\newcommand{\limit}{\textnormal{lim}}
\newcommand{\colimit}{\textnormal{colim}}

\title{Singular homology of roots of unity}

\subjclass[2020]{Primary 55N10; Secondary 54A05, 05C20} 

\author{Nikola Mili\'cevi\'c}

\begin{document}

\begin{abstract}
    We extend some basic results from the singular homology theory of topological spaces to the setting of \v{C}ech's closure spaces. We prove analogues of the excision and Mayer-Vietoris theorems and the Hurewicz theorem in dimension one. We use these results to calculate examples of singular homology groups of spaces that are not topological but are often encountered in applied topology, such as simple undirected graphs. 
    We focus on the singular homology of roots of unity with closure structures arising from considering nearest neighbors. These examples can then serve as building blocks along with our Mayer-Vietoris and excision theorems for computing the singular homology of more complex closure spaces. 
\end{abstract}

\maketitle

\section{Introduction}
A central idea in applied topology is to define and compute analogues of classical notions in algebraic topology for combinatorial objects which are typically not topological spaces~\cite{barcelo2014discrete,grigor2018path,dochtermann2009hom,lupton2022homotopy,rieser2021vcech,ebel2023synthetic}. A different approach is prevalent in topological data analysis~\cite{MR3839171,Rabadan:2019}; given an object that is not a topological space, use it to produce a filtration of topological spaces to which homology or some other functor is applied. \v{C}ech closure spaces connect these two approaches~\cite{bubenik2021eilenbergsteenrod,bubenik2024homotopy,rieser2021vcech}.

For closure spaces there exist non-constant continuous maps $f:S^1\to (\mathbb{Z}_n,c_m)$ where $c_m$ is the non-idempotent closure that expands a subset of $\mathbb{Z}_n$ to include all points within distance $m$.
This allows algebraic data structures on the sphere to be mapped on a finite sample. 

\subsection*{Our contributions}
We study certain singular homology groups of (\v{C}ech) closure spaces. A closure space is a pair $(X, c)$ where $X$ is a set and $c:2^X\to 2^X$ is a function, called a closure operator, that satisfies certain axioms (\cref{def:closure_spaces}). A singular $n$-simplex on a closure space $(X,c)$ is a continuous map $\sigma:|\Delta^n|\to (X,c)$, where $|\Delta^n|$ is the standard $n$-simplex in $\mathbb{R}^{n+1}$. The  singular chain groups are the free abelian groups generated by the singular simplices, with the obvious boundary maps. The homology of these chain groups, denoted $H_{\bullet}(X,c)$ or $H_{\bullet}(X)$ when $c$ is clear from context, is the focus of study of this manuscript.

\subsubsection*{Hurewicz theorem}
Let $I=[0,1]$ be the unit interval with the standard topology.
A homotopy relation on maps $f,g:X\to Y$ of closure spaces is defined by the existence of a map $H:X\times I\to Y$ such that $H(x,0)=f(x)$ and $H(x,1)=g(x)$. Fundamental groups can then be defined via homotopy classes of maps $f:(I,\partial I)\to (X,x_0)$ for a given basepoint $x_0$. We denote these groups by $\pi_{1}(X,x_0)$. There exists a group homomorphism, called the Hurewicz map, $h:\pi_1(X)\to H_1(X)$ sending a homotopy class of a loop to its homology class. 

\begin{theorem*}[\cref{theorem:homology_is_abelianization_of_fundamental_group}]
 If $X$ is path-connected $h$ is surjective and its kernel is the commutator subgroup of $\pi_1(X)$. Hence $H_1(X)$ is the abelianization of $\pi_1(X)$. 
\end{theorem*}

\subsubsection*{Excision and Mayer-Vietoris}
For an interior cover $\mathcal{U}$ of a closure space $X$, a $\mathcal{U}$-small singular $n$-simplex is a singular simplex, $\sigma:|\Delta^n|\to X$, whose image lies in an element of $\mathcal{U}$. We show that the chain complex of $\mathcal{U}$-small singular simplices, $C_n^{\mathcal{U}}(X)$, is chain homotopy equivalent to the chain complex of singular simplices, $C_n(X)$.

\begin{proposition*}[\cref{prop:homology_of_covering_system}]
The inclusion map $\iota:C_n^{\mathcal{U}}(X)\to C_n(X)$ is a chain homotopy equivalence. Hence $\iota$ induces isomorphisms $H_n^{\mathcal{U}}(X)\approx H_n(X)$ for all $n\in \mathbb{Z}$.
\end{proposition*}

Using this result we show that singular homology of closure spaces satisfies excision and Mayer-Vietoris theorems. 

\begin{theorem*}[\cref{theorem:excision}]
Let $\{A,B\}$ be an interior cover of the closure space $X$. For all $n\in \mathbb{Z}$, there is an isomorphism $H_n(X,A)\approx H_n(B,A\cap B)$. Equivalently, given subsets $Z\subset A\subset X$ such that the $c(Z)\subset i(A)$, the inclusion $(X\setminus Z,A\setminus Z)\to (X,A)$ induces isomorphisms $H_n(X\setminus Z,A\setminus Z)\approx H_n(X,A)$ for all $n$.
\end{theorem*}

\begin{theorem*}[\cref{theorem:Mayer_Vietoris}]
Let $\{A,B\}$ be an interior cover of the closure space $X$.
There is a long exact sequence:
\[\cdots\to H_n(A\cap B)\to H_n(A)\oplus H_n(B)\to H_n(X)\to H_{n-1}(A\cap B)\to \cdots \]
\end{theorem*}

Besides these theorems we also have results about the existence of a long exact sequence in homology for a pair $(X,A)$. For a good pair $(X,A)$ we also show that the reduced relative homology $\tilde{H}_k(X,A)$ agrees with $\tilde{H}_k(X/A)$. 

\subsubsection*{Computations}
We compute examples of singular homology for spaces that are not topological. The homology of the sphere, $S^n$, is often the first nontrivial homology computation one encounters. Using $S^n$ as a building block, the homology of more complicated spaces is calculated. This is our motivation here; to compute singular homology of simpler closure spaces that are analogous to topological spheres. Our spaces of interest are $(S^1,c_r)$ and $(\mathbb{Z}_n,c_m)$. The first one is the circle $S^1$ with the geodesic metric, $d$, and circumference $1$ yielding the closure space $(S^1,c_r)$, where 
\[c_r(A)=\{x\in S^1\,|\, \textnormal{dist}(x,A)=\inf_{y\in A}d(x,y)\le r\}, A\subset S^1,\]
for a given parameter $r\in \mathbb{R}$.
The second one, $(\mathbb{Z}_n,c_m)$, are the $n$-th roots of unity with the $m$ nearest neighbors closure operator. We have the following results.

\begin{theorem*}[\cref{theorem:homology_group_of_circle}]
$H_1(S^1,c_r)=\begin{cases}
\mathbb{Z}, & 0\le r<\frac{1}{3}\\
0, & \frac{1}{2}\le r
\end{cases}.$
\end{theorem*}

\begin{theorem*}[\cref{theorem:homology_of_roots_of_unity}]
$H_1(\mathbb{Z}_n,c_m)=\begin{cases}
\mathbb{Z}, & 3\le 3m<n\\
0, & \floor{\frac{n}{2}} \le m \text{ or } n=3, 1\le m
\end{cases}.$
\end{theorem*}

For a wedge of spaces from previous results we have the following:
\begin{theorem*}[\cref{theorem:homology_wedge_of_circles}]
For $X=(S^1,c_{r_1})\vee(S^1,c_{r_2})$,
 \[H_1(X)=\begin{cases}
\mathbb{Z}\oplus \mathbb{Z},& r_1,r_2<\dfrac{1}{3}\\
\mathbb{Z},& r_i<\dfrac{1}{3},r_j\ge \dfrac{1}{2}, i\neq j\\
0,& r_1,r_2\ge \dfrac{1}{2}
\end{cases}.\]
\end{theorem*}

\begin{theorem*}[\cref{theorem:homology_wedge_of_roots_of_unity}]
For $X=(\mathbb{Z}_{n_1},c_{m_1})\vee (\mathbb{Z}_{n_2},c_{m_2})$,
\[H_1(X)=\begin{cases}
\mathbb{Z}\oplus \mathbb{Z},& 3\le 3m_i<n_i, i=1,2 \\
\mathbb{Z},& 3m_i<n_i,\floor{\frac{n_j}{2}} \le m_j \text{ or } n_j=3, 1\le m_j, i\neq j\\
0,& \floor{\frac{n_i}{2}} \le m_i \text{ or } n_i=3, 1\le m_i, i=1,2
\end{cases}.\]
\end{theorem*}

Other constructions besides a wedge of spaces are also considered. For a closure space $X$, we show that for its suspension $\Sigma X$, $H_{k+1}(\Sigma X)=H_k(X)$. In a companion paper, it is shown that the homology of products of spaces is determined by K\"unneth theorems \cite{bubenik2021eilenbergsteenrod}.

\subsubsection*{Higher homologies}
So far we've only computed homology in degree $1$. For certain choices of parameters, homology groups in higher degrees are trivial.

\begin{theorem*}[\cref{theorem:higher_homology_for_n_greater_than_4m}]
Suppose $n\ge 4m$. Then $H_k(\mathbb{Z}_n,c_{m})=0$ for $k\ge 2$.
\end{theorem*}

\subsubsection*{Relating different roots of unity closure structures}
Finally, we give partial answers as to how the homology of roots of unity changes as we change the number of roots and the closure operation.

\begin{theorem*}[\cref{theorem:higher_homology_for_n_between_3m_and_4m_part_1}] For $m\ge 2$, suppose that $n=3m+l$, $\max(0, m-3)<l<m$. Then \[H_k(\mathbb{Z}_n,c_{m})\approx H_k(\mathbb{Z}_{n+4},c_{m+1}),\forall k\in \mathbb{N}.\]
\end{theorem*}

\begin{theorem*}[\cref{theorem:higher_homology_for_n_between_3m_and_4m_part_2}]
There exist isomorphisms
\[H_k(\mathbb{Z}_{4m-4},c_{m})\approx H_k(\mathbb{Z}_{4m},c_{m+1}),\, \forall m\ge 5,\forall k\in \mathbb{N}.\]
\end{theorem*}

\subsection*{Related work}
Extending homotopy and homology theory from topological spaces to closure spaces has been done by Demaria and Bogin \cite{demaria1985shape,demaria1984homotopy}. More recently, Rieser \cite{rieser2021vcech} developed the theory of covering spaces in the category of closure spaces and computed fundamental groups of various $(\mathbb{Z}_n,c_m)$ and $(S^1,c_r)$. The homotopy notions and homology groups studied here were also considered in \cite{bubenik2021eilenbergsteenrod,bubenik2024homotopy}.  Together with this manuscript, the results from \cite{bubenik2021eilenbergsteenrod} show that $H_n(X)$ is an Eilenberg-Steenrod homology theory for closure spaces. Palacios \cite{PalaciosVela2019} developed \v{C}ech (co)homology for closure spaces. Rieser \cite{rieser2021grothendieck} developed a novel sheaf theory for closure spaces and computed \v{C}ech and sheaf cohomologies for roots of unity $(\mathbb{Z}_n,c_1)$. There have been many approaches in applied topology that defined analogues of classical algebraic topology constructions over the last few decades. Such endeavors have been undertaken for example in the case of homology groups~\cite{arslan2008homology,barcelo2014discrete,grigor2018path,jamil2020digital,karaca2012cubical,lee2014digital,barcelo2021homology,grigor2017homologies}, homotopy theory~\cite{babson2006homotopy,barcelo2001foundations,barcelo2005perspectives,dochtermann2009hom,dochtermann2023homomorphism,kramer1998combinatorial,lupton2022homotopy,ebel2023synthetic}, fundamental groups~\cite{Boxer:1999,lupton2021fundamental,kong1989digital,lupton2022digital}, and cohomology groups~\cite{grigor2015cohomology,ege2013cohomology,gonzalez2003towards}.

\section{Background}
\label{section:background}
Here we provide background on \v{C}ech closure spaces. References are made to \cite{rieser2021vcech,vcech1966topological}. 

\subsection{Closures and interiors}
\label{section:definitions}
We start with basic definitions and facts about closure spaces.

\begin{definition}
\label{def:closure_spaces}
For a set $X$, a function $c:2^X\to 2^X$ is called a \emph{closure operator} (or \emph{closure}) for $X$ if the following axioms are satisfied:
\begin{itemize}
\item[1)] $c(\varnothing)=\varnothing$,
\item[2)] $A\subset c(A)$ for all $A\subset X$,
\item[3)] $c(A\cup B)=c(A)\cup c(B)$ for all $A,B\subset X$.
\end{itemize}
The pair $(X,c)$ where $X$ is a set and $c$ a closure for $X$ is called a (\emph{\v{C}ech}) \emph{closure space}.  A \emph{topological or Kuratowski closure operator} for a set $X$ is a closure $c$ for $X$ satisfying the additional idempotency axiom:
\begin{itemize}
\item[4)] $c(c(A))=c(A)$, for all $A\subset X$.
\end{itemize}
We say $A\subset X$ is \emph{closed} in $X$ if $c(A)=A$, and \emph{open} if $X\setminus A$ is closed.
A closure space $(X,c)$ is called a \emph{topological space} if the closure operation $c$ is \emph{topological}.
\end{definition}

\begin{remark}
\label{remark:topological_spaces_are_closure_spaces}
This definition of a topological space is equivalent to the standard definition of a topological space. More precisely, given a topology $\mathcal{O}$ on a set $X$ there is a topological closure operator $c$, induced by $\mathcal{O}$, defined by $c(A)=\{x\in X\,|\, \forall O\in \mathcal{O}, x\in O\Longrightarrow A\cap O\neq \varnothing\}$. That is $x$ is in the closure of $A$ if every open neighborhood of $x$ intersects $A$. Furthermore, if $(X,c)$ is a topological closure space, the collection of all $O$ such that $c(X\setminus O)=X\setminus O$ gives a topology on $X$ and the closure operator induced by this topology recovers the original operator $c$.
\end{remark}

There are closure spaces where the closure operation does not come from taking topological closures induced by a topology (\cref{example:closure_space_not_a_topological_space}). When $c$ is clear from context, we write $X$ for the pair $(X,c)$, however when working with multiple closure spaces at once it is  desirable to write out the closure symbols explicitly. 

\begin{example}
\label{example:closure_space_not_a_topological_space}
Consider $\mathbb{R}^n$ with euclidean metric $d$. Let $c_1$ be the closure for $\mathbb{R}^n$ defined by $c_1(A):=\{x\in \mathbb{R}^n\,|\, \textnormal{dist}(x,A):=\inf_{y\in A}d(x,y)\le 1\}$, for $\varnothing\neq A\subset\mathbb{R}^n$. Let $U=\{x\in \mathbb{R}^n\,|\,d(x,0)\le 1\}$ be the radius $1$ ball around the origin. Then $c_1(U)$ is a radius $2$ ball around the origin, while the double closure $c_1(c_1(U))$ is a radius $3$ ball around the origin. Therefore $c_1^2\neq c_1$, and hence $c_1$ is not a topological closure.
\end{example}

Closure operations are monotone; it is elementary that  $A\subset B\subset X \Longrightarrow c(A)\subset c(B)$. 

\begin{example}
\label{example:discrete_and_indiscrete_closures}
Let $X$ be a set. The identity $\id{2^X}:2^X\to 2^X$ is a topological closure for $X$. It is the discrete closure for $X$ and it is the closure of the discrete topology on $X$. The topological closure operation defined by $A\mapsto X$ for $A\neq \varnothing$ and $\varnothing\mapsto\varnothing$ is called \emph{the indiscrete closure} for $X$ and it corresponds to the closure of the indiscrete topology we can put on $X$.
\end{example}

Attached to a closure operation comes an interior operation.

\begin{definition}
\label{def:interior}
For a closure space $(X,c)$ there is an associated \emph{interior operation}, $i_c:2^X\to 2^X$ given by 
\[i_c(A):=X\setminus c(X\setminus A),\, \forall A\subset X.\]
The interior operation satisfies the following:
\begin{itemize}
\item[1)] $i_c(X)=X$,
\item[2)] For all $A\subset X$, $i_c(A)\subset A$,
\item[3)] For all $A,B\subset X$, $i_c(A\cap B)=i_c(A)\cap i_c(B)$.
\end{itemize}
If $i:2^X\to 2^X$ is a function satisfying the 3 conditions above and we define $c_i:2^X\to 2^X$ by 
\[c_i(A):=X\setminus i(X\setminus A),\, \forall A\subset X,\]
then $c_i$ is a closure operation for $X$. When $c$ is clear from context we will write $i$ for $i_c$.

\begin{proposition}{\cite[14.A.12]{vcech1966topological}}
\label{proposition:open_if_equal_to_interior}
$A\subset X$ is open in $X$ if and only if $i(A)=A$. \qed
\end{proposition}
\end{definition}

\subsection{Continuous functions}
Here we recall properties of continuous maps of closure spaces. For a more detailed background see~\cite{vcech1966topological}.

Let $(X,c)$ and $(Y,d)$ be closure spaces. A set map $f:(X,c)\to (Y,d)$ is \emph{continuous} if for every $A\subset X$, $f(c(A))\subset d(f(A))$. A continuous map $f$ is called a \emph{homeomorphism} if $f$ is a bijection and $f^{-1}$ is continuous. We have that  $f(c(A))=d(f(A))$ for all $A\subset X$ if and only if $f$ is a homeomorphism. Compositions of continuous maps are continuous.

For $A\subset X$, the \emph{subspace closure operation} on $A$ is $c_A(B):=A\cap c(B)$ for all $B\subset A$. We say that $(A,c_A)$ is a \emph{subspace} of $(X,c)$. The inclusion map $(A,c_A)\to (X,c)$ is continuous. Restrictions of continuous maps to subspaces are continuous. Let $A\subset B\subset X$. We say that $B$ is a \emph{neighborhood} of $A$ in $X$ if $A\subset i(B)$. If $A=\{x\}$ is a singleton, we say $B$ is a neighborhood of $x$.

Let $c_1$ and $c_2$ be two closure operations for $X$. We say $c_1$ is \emph{coarser} than $c_2$ and $c_2$ is \emph{finer} than $c_1$ if $c_2(A)\subset c_1(A)$ for all $A\subset X$. Equivalently, as $(2^X,\subset)$ is a poset, $c_1$ is coarser than $c_2$ if and only if $c_1$ dominates $c_2$ as a poset morphism. Observe  that $c_1$ is coarser than $c_2$ iff the identity map $\id{X}:(X,c_2)\to (X,c_1)$ is continuous.

We denote the category with objects closure spaces and morphisms continuous maps of closure spaces by $\cat{Cl}$ throughout. We denote the full subcategory of $\cat{Cl}$ with objects topological spaces by $\cat{Top}$. There is a pair of adjoint functors between $\cat{Cl}$ and $\cat{Top}$.

\begin{proposition}{\cite[16.B.1-16.B.3]{vcech1966topological}}
\label{prop:topological_modification}
Let $(X,c)$ be in $\cat{Cl}$, let $\tau (c):2^X\to 2^X$ be defined by 
\[\tau(c)(A):=\bigcap\{F\subset X\,|\, c(F)=F\text{ and } A\subset F\}.\]
Then $\tau (c)$ is a topological closure operation, called \emph{the topological modification of $c$}. Moreover, $\tau (c)$ is the finest topological closure operation coarser than $c$. \qed
\end{proposition}

\begin{proposition}{\cite[16.B.4]{vcech1966topological}}
\label{prop:adjoints_between_closures_and_topologies}
Let $(X,c)$ be a closure space and let $(Y,d)$ be a topological space. A set theoretic map $f:X\to Y$ is continuous as a closure space map $f:(X,c)\to (Y,d)$ if and only if the map $f:(X,\tau(c))\to (Y,d)$ is a continuous map of topological spaces. That is to say there exists a natural bijection between the sets of morphisms
\[\cat{Cl}((X,c),(Y,d))\approx \cat{Top}((X,\tau(c)),(Y,d)).\]
Equivalently, $\tau$ is the left adjoint to the inclusion functor $\iota:\cat{Top}\to \cat{Cl}$. \qed
\end{proposition}

\subsection{Covers of closure spaces} 
Here we recall different types of covers of closure spaces. Let $X$ be a closure space.

A family $\{U_{\alpha}\, |\, \alpha\in A \}$ of subsets of $X$ is called \emph{locally finite} if each point $x\in X$ possesses a neighborhood intersecting only finitely many $U_{\alpha}$. A \emph{cover} is a family of subsets of $X$, $\mathcal{U}=\{U_{j}\}_{j\in J}$, whose union is $X$. $\mathcal{U}$ is an \emph{interior cover} of $X$ if every $x\in X$ has a neighborhood in $\mathcal{U}$. Equivalently, $\mathcal{U}$ is an interior cover of $X$ if $X = \bigcup_{U\in \mathcal{U}} i(U)$, where $i$ is the interior operator for $X$.
We say $\mathcal{U}$ is an \emph{open (closed) cover} if every $U_j$ is open (closed). If $X$ is a topological space, any open cover of $(X,c)$ is an interior cover. We have the following version of the pasting lemma.

\begin{theorem}{\cite[17.A.18]{vcech1966topological}}[Pasting Lemma]
\label{theorem:pasting_lemma_for_closure_spaces}
Let $(X,c)$ be a closure space and let $\{U_{\alpha}\,|\,\alpha\in A\}$ be a locally finite closed cover of $X$. Let $f:(X,c)\to (Y,d)$ be a set-theoretic map. If $f|_{U_{\alpha}}$ is continuous for each $\alpha\in A$, then $f$ is continuous. \qed
\end{theorem}

\begin{lemma}
\label{lemma:finite_cover_whose_image_is_contained_in_covering_system}
Let $(X,c)$ be a compact topological space. Let $(Y,d)$ be a closure space with interior cover $\mathcal{V}$, and suppose that $f:(X,c)\to (Y,d)$ is continuous. Then there exists an (finite) open cover $\mathcal{U}$ of $(X,c)$ such that for all $U\in \mathcal{U}$, there is a $V\in \mathcal{V}$ such that $f(U)\subset V$.
\end{lemma}

To prove \cref{lemma:finite_cover_whose_image_is_contained_in_covering_system} we will be relying on the following result. 

\begin{theorem}{\cite[Theorem 16.A.4 and Corollary 16.A.5]{vcech1966topological}}
\label{theorem:equivalent_definitions_of_continuity}
Let $(X,c)$ and $(Y,d)$ be closure spaces. A map $f:(X,c)\to (Y,d)$ is continuous at $x\in X$ if and only if for every neighborhood $V\subset Y$ of $f(x)$, the inverse image $f^{-1}(V)$ is a neighborhood of $x$. Equivalently, $f$ is continuous at $x$ if and only if for each neighborhood $V\subset Y$ of $f(x)$, there exists a neighborhood $U\subset X$ of $x$ such that $f(U)\subset V$. \qed
\end{theorem}

\begin{proof}[Proof of \cref{lemma:finite_cover_whose_image_is_contained_in_covering_system}]
    Since $\mathcal{V}$ is an interior cover, for all $x\in X$, there is a neighborhood of $f(x)$ in $\mathcal{V}$. Label this neighborhood by $V_{f(x)}$. By \cref{theorem:equivalent_definitions_of_continuity}, the neighborhood $V_{f(x)}\subset Y$ of $f(x)$ there exists a neighborhood $W_x \subset X$ of $x$ in $(X,c)$ such that $f(W_x) \subset V_{f(x)}$. Since $W_x$ is a neighborhood of $x$ and $X$ is topological, there is an open set $U_x \subset W_x$ with $x \in U_x$ for every $x \in X$. The collection $\{U_x\}_{x\in X}$ is thus an open cover of $X$ which refines $\{W_x\}_{x\in X}$ by construction. Since $X$ is compact, $\{U_x\}_{x\in X}$ admits a finite subcover $\mathcal{U}=\{U_{x_i}\}_{i=1}^n$, and $\mathcal{U}$ satisfies the conclusion of the lemma, by construction. 
\end{proof}

\subsection{Limits and colimits of closure spaces}
\label{section:canonical_closure_operations}

Here we define limits and colimits of diagrams of closure spaces. Then we construct products, coproducts, equalizers and coequalizers of closure spaces which implies that $\cat{Cl}$ is complete and cocomplete. For a more detailed background on limits and colimits in $\mathbf{Cl}$ see~\cite{vcech1966topological}.

Let $X:\mathcal{I}\to \cat{Cl}$ be a diagram of closure spaces, i.e. $X$ is a functor from some indexing category $\mathcal{I}$ to $\cat{Cl}$. In particular, for all objects $i$ in $\mathcal{I}$ we have a closure space $X_i$ and for every morphism $\varphi:i\to j$ in $\mathcal{I}$ we have a continuous map $X_{\varphi}:X_i\to X_j$.

\begin{definition}
    A \emph{limit} of the $\mathcal{I}$-diagram $X$ in $\cat{Cl}$ is given by a closure space $\limit_{\mathcal{I}}X$ together with continuous maps $p_i:\lim_{\mathcal{I}}X\to X_i$ such that
    \begin{enumerate}[left=0pt]
        \item For $\varphi:i\to j$ a morphism in $\mathcal{I}$ we have $p_j=X_{\varphi}\circ p_i$.
        \item For any closure space $Y$ and any family of continuous maps $q_i:Y\to X_i$, indexed by $\mathcal{I}$, such that for all $\varphi:i\to j$ in $\mathcal{I}$ we have $q_j=X_{\varphi}\circ q_i$ there exists a unique morphism $q:Y\to \lim_{\mathcal{I}}X$ such that $q_i=p_i\circ q$ for every object $i$ in $\mathcal{I}$.
    \end{enumerate}
\end{definition}

\begin{example}
    Let $I$ be a set. Let $\mathcal{I}$ be the category whose objects are the elements of $I$ and whose morphisms are only the identity morphisms. For a diagram of closure spaces $X:\mathcal{I}\to \cat{Cl}$, $\limit_{\mathcal{I}}X_i$ is the product $\prod_{i\in I}X_i$ with projections maps $p_i:\prod_{i\in I}X_i\to X_i$ and the closure on the underlying set of $\prod_{i\in I}X_i$ is the coarsest closure making all the projection maps continuous. We can also describe this closure in terms of neighborhoods in the following way.
    For all $x\in \prod_{i\in I}X_i$, let $\mathcal{U}_x$ denote the collection of sets of the form
    \begin{equation*} \label{eq:product}
        \bigcap\{p_i^{-1}(V_i)\,|\,i\in A\},
    \end{equation*}
    where $A\subset I$ is finite and $V_i$ is a neighborhood of $p_i(x)$ in $X_i$. 
    For $B\subset  X$, let 
    \[c(B)=\{x\in X\,|\, \forall U\in \mathcal{U}_x, U\cap B\neq \varnothing\}.\]
    The closure $c$ is called the \emph{product closure} for $\prod_{i\in I}X_i$ and $(\prod_{i\in I}X_i,c)$ is called the \emph{product closure space}. 
    For closure spaces $(X,c)$ and $(Y,d)$, their product is denoted by $(X\times Y,c\times d)$.
\end{example}

\begin{example}
    Let $\mathcal{I}$ be the category with two objects, $1$ and $2$ and two parallel morphisms from $1$ to $2$. Then, an $\mathcal{I}$-diagram of closure spaces $X:\mathcal{I}\to \cat{Cl}$ consists of a pair of continuous maps
    $
    \begin{tikzcd}(X_1,c_1) \ar[r,shift left,"f"] \ar[r,shift right,"g"'] & (X_2,c_2).
    \end{tikzcd}$
    We call $\limit_{\mathcal{I}}X$ the \emph{equalizer of $f$ and $g$}. In particular, the equalizer consists of the closure space $(E,c_E)$ and map $\iota:E \to X_1$ defined in the following manner. Let $E=\{x \in X_1 \ | \ f(x)=g(x)\}$ with $\iota$ the inclusion map. For any $A \subset E$, set $c_E(A) = c(A) \cap E$.
\end{example}

We thus have the following.

\begin{theorem}{\cite[Theorems 32.A.4 and 32.A.10]{vcech1966topological}} \label{complete}
  The product and equalizer defined above are the categorical product and equalizer in $\cat{Cl}$ and hence $\cat{Cl}$ is complete. \qed
\end{theorem}

Dually, we also define colimits of diagrams closure spaces.

\begin{definition}
    A \emph{colimit} of the $\mathcal{I}$-diagram $X$ in $\cat{Cl}$ is given by a closure space $\colimit_{\mathcal{I}}X$ together with continuous maps $s_i:X_i\to \colimit_{\mathcal{I}}X$ such that
    \begin{enumerate}[left=0pt]
        \item For $\varphi:i\to j$ a morphism in $\mathcal{I}$ we have $s_i=s_j\circ X_{\varphi}$.
        \item For any closure space $Y$ and any family of continuous maps $t_i:X_i\to Y$, indexed by $\mathcal{I}$, such that for all $\varphi:i\to j$ in $\mathcal{I}$ we have $t_i=t_j\circ X_{\varphi}$ there exists a unique morphism $t:\colimit_{\mathcal{I}}X\to Y$ such that $t_i=t\circ s_i$ for every object $i$ in $\mathcal{I}$.
    \end{enumerate}
\end{definition}

\begin{example}
    Let $I$ be a set. Let $\mathcal{I}$ be the category whose objects are the elements of $I$ and whose morphisms are only the identity morphisms. For a diagram of closure spaces $X:\mathcal{I}\to \cat{Cl}$, $\colimit_{\mathcal{I}}X_i$ is the coproduct $\coprod_{i\in I}X_i$ with coprojections maps $s_i:X_i\to \coprod_{i\in I}X_i$ and the closure on the underlying set of $\coprod_{i\in I}X_i$ is the finest closure making all the coprojection maps continuous. More specifically, the \emph{coproduct}
    of $\{(X_i,c_i)\}_{i\in I}$ is the  disjoint union of sets $X=\sqcup_i X_i$ with the closure operation $c$ for $X$ defined by $c(\sqcup_i A_i):=\sqcup_i c_i(A_i)$ for all subsets $\sqcup_i A_i$ of $\sqcup_{i\in I}X_i$.
\end{example}

\begin{example}
\label{ex:coequalizer}
    Let $\mathcal{I}$ be the category with two objects, $1$ and $2$ and two parallel morphisms from $1$ to $2$. Then, an $\mathcal{I}$-diagram of closure spaces $X:\mathcal{I}\to \cat{Cl}$ consists of a pair of continuous maps $
    \begin{tikzcd}(X_1,c_1) \ar[r,shift left,"f"] \ar[r,shift right,"g"'] & (X_2,c_2).
    \end{tikzcd}
    $
    We call $\colimit_{I}X$ the \emph{coequalizer of $f$ and $g$}. In particular, the coequalizer consists of the closure space $(Q,c_q)$ and a map $q:X_2 \to Q$ defined in the following manner. Let $Q$ be the quotient set $X_2/\!\! \sim$, where $\sim$ is the smallest equivalence relation such that  $f(x) \sim g(x)$, for all $x \in X_1$. Let $q:X_2 \to Q$ be the quotient map. For any $A \subset Q$, set $c_q(A) = q(c_2(q^{-1}(A)))$.
\end{example}

\begin{example}
\label{ex:pushout}
    Let $\mathcal{I}$ be the category with $3$ objects, $0$, $1$ and $2$ and the only non-identity morphisms being $0\to 1$, and $0\to 2$. Let $X$ be an $\mathcal{I}$-diagram of closure spaces, $X:\mathcal{I}\to \cat{Cl}$. Then we have non-identity morphisms $f:(X_0,c_0)\to (X_1,c_1)$ and $g:(X_0,c_0)\to (X_2,c_2)$, in the diagram.
    We call $\colimit_{\mathcal{I}}X$ the \emph{pushout} of $f$ and $g$. More specifically, the pushout is given by the closure space $(Y,d)$, along with induced maps $s_1:(X_1,c_1)\to (Y,d)$ and $s_2:(X_2,c_2)\to (Y,d)$ in the following commutative diagram:
    \begin{center}
        \begin{tikzcd}
        (X_0,c_0)\arrow[r,"f"]\arrow[d,"g"] & (X_1,c_1)\arrow[d,"s_1",dashed]\\
        (X_2,c_2)\arrow[r,"s_2",dashed] & (Y,d)
    \end{tikzcd}
    \end{center}
    where we set $Y=(X_1\sqcup X_2)/\sim$, with $\sim$ being the smallest equivalence relation such that for all $x\in X_0$, $f(x)\sim g(x)$, and the closure operator $d$ defined by $d(A)=s_1(c_1(s_1^{-1}(A)))\cup s_2(c_2(s_2^{-1}(A)))$, for all $A\subset Y$.
\end{example}

\begin{theorem}{\cite[Theorems 33.A.4 and 33.A.5]{vcech1966topological}} \label{cocomplete}
 The coproduct and coequalizer defined above are the categorical coproduct and coequalizer in $\cat{Cl}$ and hence $\cat{Cl}$ is cocomplete. \qed
\end{theorem}

\begin{proposition}
\label{prop:limits_and_colimits_of_closure_spaces}
Every limit of topological spaces in $\cat{Cl}$ is a topological space. On the other hand, colimits of topological spaces in $\cat{Cl}$ are not necessarily topological spaces.
\end{proposition}

\begin{proof}
The inclusion functor $\iota:\cat{Top}\to \cat{Cl}$ is a right-adjoint by \cref{prop:adjoints_between_closures_and_topologies} and thus preserves limits. For the second claim, see \cite[Introduction to Section 33.B]{vcech1966topological}.
\end{proof}

\begin{remark}
\label{remark:adjunction_top_and_cl}

If $X$ is a colimit of a diagram of topological spaces in $\cat{Top}$, it may be a closure space and not a topological space in $\cat{Cl}$, as $\cat{Top}$ is a reflective, but not coreflective, subcategory of $\cat{Cl}$ (\cref{prop:limits_and_colimits_of_closure_spaces}). To differentiate the two, let $X_{\cat{Top}}$ and $X_{\cat{Cl}}$ be the colimits in $\cat{Top}$ and $\cat{Cl}$, respectively, of a given diagram of topological spaces. There is a continuous map induced by the unit of the adjunction (\cref{prop:adjoints_between_closures_and_topologies}) and it is the set theoretic identity map $\id{X}:X_{\cat{Cl}}\to X_{\cat{Top}}$. That is, $X_{\cat{Top}}$ is the topological modification of $X_{\cat{Cl}}$. In other words, $X_{\cat{Cl}}$ and $X_{\cat{Top}}$ are closure spaces on the same underlying set, and $X_{\cat{Top}}$ is the finest topological space coarser than $X_{\cat{Cl}}$. 
\end{remark}

For some diagrams of topological spaces of interest, $X_{\cat{Cl}}$ and $X_{\cat{Top}}$ from \cref{remark:adjunction_top_and_cl} coincide, see \cref{theorem:cw_complexes} below.

\subsection{CW complexes in $\cat{Cl}$}
Here we briefly recall how CW complexes are constructed in $\cat{Cl}$. For a more detailed study of CW complexes in $\cat{Cl}$, see \cite{bubenik2023cw}.

Let $S^n$ be the topological $n$-sphere and $D^n$ the topological $n$-disk. A closure space $(X,c)$ is obtained from a closure space $(A,b)$ by \emph{attaching cells} if there exists a pushout in $\cat{Cl}$,

\begin{center}
    \begin{tikzcd}
        \coprod_{j} S_j^{n_j-1}\arrow[d,hookrightarrow,"\coprod_j\iota_j"]\arrow[r,"\varphi"] & (A,b)\arrow[d,"\iota_A",dashed]\\
        \coprod_{j}D_j^{n_j}\arrow[r,"\Phi",dashed] & (X,c)
    \end{tikzcd}
\end{center}
where for each $j$, $D_j^{n_j}$ is the $n_j$-dimensional disk and $\iota_j:S_j^{n_j-1}\to D_j^{n_j}$ is the inclusion of its boundary. We call $\varphi$ the \emph{cell-attaching map}. If for all $j$, $n_j=n$ we say that $(X,c)$ is obtained from $(A,b)$ by \emph{attaching $n$-cells}.

For a closure space $X$ and a closed subspace $A$, we say the pair $(X,A)$ is a \emph{relative CW complex} in $\cat{Cl}$ if $X$ is the colimit of a diagram of closure spaces
\[A=X^{-1}\hookrightarrow X^0\hookrightarrow X^1\hookrightarrow X^2\hookrightarrow\cdots,\]
where for all $n\ge 0$, $X^n$ is obtained from $X^{n-1}$ by attaching $n$-cells. If $A=\varnothing$, we say $X$ is a \emph{CW complex} in $\cat{Cl}$. If the total numbers of cells attached is finite, we call $X$ a \emph{finite CW complex}. If for each $n\ge 0$, there are only finitely many $n$-cells attached, we call $X$ a \emph{CW complex of finite type}. If no more cells are attached after a certain $n\ge 0$, we say $X$ is a \emph{finite-dimensional CW complex.}

\begin{theorem}{\cite[Theorem 1.4]{bubenik2023cw}}
\label{theorem:cw_complexes}
    The constructions of finite CW complexes relative to a compactly generated weak Hausdorff topological space in $\cat{Cl}$ and $\cat{Top}$ agree. \qed
\end{theorem}

The same is not true for finite-dimensional CW complexes, CW complexes of finite type, and finite relative CW complexes, see \cite{bubenik2023cw} for examples. 

\subsection{Homotopy for closure spaces}
\label{section:continuous_theories}
Here we recall a homotopy theory for closure spaces and the definition of higher homotopy groups. These have been considered in \cite{demaria1985shape,demaria1984homotopy,rieser2021vcech}. In \cite{bubenik2021eilenbergsteenrod,bubenik2024homotopy} several different homotopy theories for closure spaces were defined, and this particular homotopy theory was called $(I_{\tau},\times)$-homotopy, $I_{\tau}$ being the unit interval $[0,1]$ with the standard topology. We simplify this notation in this manuscript and from now on, unless otherwise specified, $I$ will denote the unit interval with the standard topology.

\subsubsection{Homotopy} 
\label{section:homotopy}
Here we recall notions homotopy and deformation retractions for closure spaces. See~\cite{bubenik2021eilenbergsteenrod,rieser2021vcech} for more details. We show cones are contractible.

Let $f,g:X\to Y$ be maps of closure spaces. We say $f$ and $g$ are \emph{homotopic} and write $f\simeq g$, if there exists a continuous function $H:X\times I\to Y$ such that $H(x,0)=f(x)$ and $H(x,1)=g(x)$. The map $H$ is called a \emph{homotopy} between $f$ and $g$. We say two closure spaces $X$ and $Y$ are homotopy equivalent if there exist continuous maps $f:X\to Y$ and $g:Y\to X$ such that $fg\simeq \id{Y}$ and $gf\simeq \id{X}$.

A closure space $X$ is \emph{path-connected} if for any $x,y\in X$ there exists a continuous map $\alpha:I\to X$ such that $\alpha(0)=x$ and $\alpha(1)=y$. We say $\alpha$ is a path between $x$ and $y$. Equivalently a path between $x$ and $y$ is a homotopy between the constant maps $*\to X$ whose images are $x$ and $y$, respectively. We put a relation on $X$ by saying $x\sim y$ if there is a continuous $\alpha:I\to X$ such that $\alpha(0)=x$ and $\alpha(1)=y$. By reversing and concatenating paths (made possible by \cref{theorem:pasting_lemma_for_closure_spaces}) one can show $\sim$ is an equivalence relation on $X$. We call the equivalence classes of $\sim$ the \emph{path components} of $X$.

Let $X$ be a closure space and let $A$ be a subspace of $X$. A map $r:X\to X$ is a \emph{retraction} of $X$ onto $A$ if $r(X)=A$ and $r|_A=\id{A}$, and we say $A$ is a retract of $X$. A \emph{deformation retraction} is a retraction $r:X\to X$ along with a homotopy between the identity and $r$. In this case, we say $X$ \emph{deformation retracts onto $A$}. $X$ is \emph{contractible} if $X$ deformation retracts to a one point space. 

\begin{lemma}
\label{lemma:retracts_of_contractible_are_contractible}
Let $X$ be a closure space and $A\subset X$ a subspace that is also a retract of $X$. If $X$ is contractible, then so is $A$.
\end{lemma}

\begin{proof}
Let $r:X\to X$ be a retraction of $X$ onto $A$. Let $f:X\to \{*\}$ and $g:\{*\}\to X$ be continuous maps in $\cat{Cl}$ such that $fg=\textnormal{id}_{\{*\}}$ and $gf\simeq \textnormal{id}_{X}$. Let $f'=f|_{A}$ be the restriction of $f$ to $A$ and let $g'=rg$. As there is only one map from $\{*\}$ to $\{*\}$ we have $f'g'=\textnormal{id}_{\{*\}}$. All that is left to show is $g'f'\simeq \textnormal{id}_{A}$. Let $H:X\times I\to X$ be a homotopy between $gf$ and $\textnormal{id}_X$. Let $H'=rH|_{A\times I}$. Then $H'$ is the desired homotopy between $g'f'$ and $\textnormal{id}_A$. 
\end{proof}

Call $(X,A)$ a \emph{good pair} if there is a neighborhood $B$ of $A$ in $X$ that deformation retracts onto $A$.

\begin{proposition}
\label{prop:maps_from_quotient_closure}
Let $(X,c_X)$ be a closure space and let $E$ be an equivalence relation on $X$ and let $q$ be the canonical quotient map, $q:(X,c_X)\to (X/E,c_q)$. A mapping $f:(X/E,c_q)\to (Y,c_Y)$ is continuous if and only if $f\circ q:(X,c_X)\to (Y,c_Y)$ is continuous. 
\end{proposition}

\begin{proof}
    Consider $E\subset X\times X$ with the induced closure, as a subspace of the closure space $(X\times X,c_X\times c_X)$. Let $p_1,p_2:E\to (X,c_X)$ be the canonical projections. Then, the quotient mapping $q:(X,c_X)\to (X/E,c_q)$ is the coequalizer of the diagram
    \begin{tikzcd}
        E \ar[r,shift left,"p_1"] \ar[r,shift right,"p_2"'] & (X,c_X),
    \end{tikzcd} and from this observation the result immediately follows.
\end{proof}

For a closure space $(X,c)$, let $(CX,c_q)$ be the cone of $(X,c)$, namely $CX=X\times I/\sim$ where $(x,0)\sim (y,0)$ for all $x,y\in X$ and $X\times I$ has the product closure operation $c\times \tau$ while by $q$ we denote the quotient map and by $c_q$ the quotient closure operation.
Cones over a space are, as expected, contractible.

\begin{proposition}
\label{prop:cone}
For a closure space $(X,c)$, $(CX,c_q)$ is contractible.
\end{proposition}

\begin{proof}
Define $H:(X\times I\times I,c\times \tau\times \tau)\to (X\times I,c\times \tau)$ by $H(x,t,s):=(x,(1-s)t+s)$. Note that $H$ is continuous. Define $\tilde{H}:(CX\times I,c_q\times \tau)\to (CX,c_q)$ by $\tilde{H}([x,t],s):=[x,(1-s)t+s]$. Note that by construction $q\circ H=\tilde{H}\circ q$. Furthermore $q\circ H$ is continuous, thus $\tilde{H}\circ q$ is continuous. By \cref{prop:maps_from_quotient_closure}, the map $\tilde{H}$ is continuous. Now observe that by construction, $\tilde{H}$ is a deformation retraction of $CX$ to a one point space, thus the result follows.
\end{proof}

\subsubsection{Homotopy extension property.}
Let $X$ be a closure space and $A\subset X$ a subspace. We say the pair $(X,A)$ has the \emph{homotopy extension property} if every pair of maps $X\times \{0\}\to Y$ and $A\times I\to Y$ that agree on $A\times \{0\}$ can be extended to a map $X\times I\to Y$. 

\begin{lemma}{\cite[Lemma 4.13]{rieser2021vcech}}
\label{lemma:hep}
If a pair of topological Hausdorff spaces $(X,A)$ satisfies the homotopy extension property in $\cat{Top}$, then $(X,A)$ also satisfies the homotopy extension property in $\cat{Cl}$. \qed
\end{lemma}

\begin{proposition}
\label{prop:quotient_map_is_homotopy_equivalence}
If the closure space pair $(X,A)$ satisfies the homotopy extension property and $A$ is contractible, then the quotient map $q:X\to X/A$ is a homotopy equivalence.
\end{proposition}

\begin{proof}
The same arguments as in the proof of \cite[Proposition 0.17]{MR1867354} apply.
\end{proof}

\subsubsection{The fundamental group} 
Finally, we recall the notion of continuous maps between pairs of closure spaces and how homotopy classes of such maps can be used to define fundamental groups of closure spaces. One could also in a straightforward way generalize this construction and define higher homotopy groups of closure spaces, but we choose not to pursue this direction here. References are made to~\cite{rieser2021vcech,demaria1984homotopy}.

A \emph{continuous map} $f:(X,A)\to (Y,B)$ between closure space pairs is a continuous map $f:X\to Y$ such that $f(A)\subset B$. For a point $p$ in $X$, we call the pair $(X,p)$ a \emph{pointed space.}

Let $X$ be a closure space and let $p\in X$, and let $I$ be the unit interval, with $\partial I$ being its boundary, the set $\{0,1\}$. The \emph{fundamental group} of the pointed space $(X,p)$ is the set 
\[\pi_1(X,p)=[(I,\partial I),(X,p)],\] 
of homotopy classes of maps $f:I\to X$ such that $f(\partial I)=\{p\}$, with the group structure described below. Let $[f], [g]\in \pi_{1}(X,p)$. We let $[f]\cdot[g]$ be represented by the homotopy class of $f\cdot g:I\to X$, where
\[f\cdot g(t)=\begin{cases}
    f(2), & t\le \frac{1}{2},\\
    g(2t-1), & \frac{1}{2}\le t.
\end{cases}\]
This operation is indeed a group operation, and this whole construction is an extension of the usual construction of fundamental groups for topological spaces, one can find for example in \cite[Chapter 6]{tom2008algebraic}. Finally, if $X$ is a path-connected closure space, then $\pi_1(X,p)$ will not depend on the basepoint $p$, and we can simply write $\pi_1(X)$ instead. Continuous maps $f:(I,\partial I)\to (X,p)$ are also called \emph{loops} in $X$ based at $p$, and $f\cdot g$, is called the \emph{concatenation} of the loops $f$ and $g$. Given a path (resp. loop) $f:I\to X$, its \emph{inverse path} (resp. \emph{inverse loop}) is defined by $\overline{f}(t)=f(1-t)$. Just like in the case of topological spaces, for a loop $f$ based at $p$, $[f]$ and $[\overline{f}]$ are inverses in $\pi_1(X,p)$.

Let $p_1,p_2\in X$ be in the same path-component and let $g:I\to X$ be a path between $p_1$ and $p_2$. For each loop $f$ based at $p_2$, we can associate the loop $g\cdot f\cdot \overline{g}$ based at $p_1$ (either the loop $(g\cdot f)\cdot \overline{g}$ or the loop $g\cdot(f\cdot\overline{g})$ but one can show they are in the same homotopy class). Let $\beta_g:\pi_1(X,p_2)\to \pi_1(X,p_1)$ be defined by $\beta_g([f])=[g\cdot f\cdot \overline{f}]$ which is well defined as it sends homotopies of loops to homotopies of loops.

\begin{lemma}
    \label{lemma:fundamental_group_path_connected}
    The map $\beta_g:\pi_1(X,x_2)\to \pi_1(X,x_1)$ is a group isomorphism.
\end{lemma}

\begin{proof}
    The function $\beta_g$ is a group homomorphism since $\beta_g([f_1\cdot f_2])=[g\cdot f_1\cdot f_2\cdot \overline{g}]=[g\cdot f_1\cdot \overline{g}\cdot g\cdot f_2\cdot \overline{g}]=\beta_g([f_1])\beta_g([f_2])$. Furthermore, $\beta_g$ has an inverse $\beta_{\overline{g}}$ since $\beta_g\beta_{\overline{g}}([f])=\beta_g([\overline{g}\cdot f\cdot g])=[g\cdot \overline{g}\cdot f\cdot g\cdot \overline{g}]=[f]$ and similarly $\beta_{\overline{g}}\beta_g([f])=[f]$.
\end{proof}

Therefore, for path connected spaces the group $\pi_1(X,p)$ does not depend on the basepoint $p$ and we abbreviate  $\pi_1(X,p)$ to $\pi_1(X)$.

\subsection{Singular homology for closure spaces}
\label{section:homology}
The $n$-\emph{dimensional standard simplex} is 
\[|\Delta^n|=\{(t_0,\dots, t_n)\in \mathbb{R}^{n+1}\,|\,\sum_it_i=1,t_i\ge 0\}.\]
Alternatively, $|\Delta^n|$ is the convex hull of the $e_i$'s, the standard basis vectors of $\mathbb{R}^{n+1}$, which we write as $[e_0,\dots, e_n]$.
Let $X$ be a closure space. Define
\[C_n(X):=\mathbb{Z}\langle\{\sigma:|\Delta^n|\to X\}\,|\, \sigma\,\text{is continuous}\rangle,\]
where $|\Delta^n|$ inherits the subspace topology from $\mathbb{R}^{n+1}$ with the Euclidean topology. The boundary maps $\partial_n:C_n(X)\to C_{n-1}(X)$ are defined by 
\[\partial_n(\sigma)=\sum_{i=1}^n(-1)^i\sigma|_{[e_0,\dots ,e_{i-1},\hat{e}_i,e_{i+1},\dots, e_n]}\]
 and extended linearly. We call $\sigma:|\Delta^n|\to X$ a \emph{singular} $n$-simplex, and elements of $C_n(X)$ \emph{singular} $n$-chains. Write $H_n(X)$ for the homology groups of this chain complex. That is $H_n(X)=\ker \partial_n/\textnormal{im }\partial_{n+1}$. We also call elements of $\ker\partial_n$ the $n$-cycles.

This definition is an extension of singular homology of topological spaces. Here we observe some properties these homology groups have. In \cite{bubenik2021eilenbergsteenrod,bubenik2024homotopy} these homology groups were denoted by $H^{I_{\tau}}_{n}(X)$. Additionally this extension of singular homology to closure spaces was first considered in \cite{demaria1984homotopy}, to the best of our knowledge.

From the definition it follows that the homology of a space can be recovered from the homology of its path components, and the $0$-th homology counts the number of path components.

\begin{proposition}
\label{prop:decomposition_into_path_components}
Let $X=\bigcup_{\alpha\in A}X_i$ be a closure space with a decomposition into path-components. Then, for all $n\in \mathbb{N}$ there is an isomorphism $H_n(X)\approx \bigoplus\limits_{\alpha\in A} H_n(X_i)$. If $X\neq\varnothing$ in $\cat{Cl}$ is path-connected, then $H_0(X)\approx \mathbb{Z}$. \hfill \qed
\end{proposition}

Let $X\neq\varnothing$ be a closure space. We define the \emph{reduced homology groups} $\tilde{H}_n(X)$ to be the homology groups of the augmented chain complex
\[\cdots\to C_2(X)\xrightarrow{\partial_2} C_1(X)\xrightarrow{\partial_1} C_0(X)\xrightarrow{\epsilon} \mathbb{Z}\to 0\]
where $\epsilon(\sum_i n_i\sigma_i)=\sum_i n_i$. It follows that $H_0(X)\approx \tilde{H}_0(X)\oplus \mathbb{Z}$ and  $H_n(X)\approx \tilde{H}_n(X)$, for all $n\ge 1$. 

Let $f:X\to Y$ be a map of closure spaces. If $\sigma:|\Delta^n|\to X$ is a singular simplex on $X$, then $f\circ\sigma$ is a singular simplex on $Y$. Thus we have an induced chain map $f_{\#}:C_n(X)\to C_n(Y)$ and in turn an induced map $f_*:H_n(X)\to H_n(Y)$ on homology.  

\begin{theorem}
\label{theorem:homotopic_maps_induce_same_homomorphism}
If  $f,g:X\to Y$ are homotopic, then they induce the same homomorphisms $f_*=g_*:H_n(X)\to H_n(Y)$ for all $n\in \mathbb{N}$.
\end{theorem}

\begin{proof}
    The proof is a verbatim replication of the standard one used for topological spaces, see for example \cite[Theorem 9.3.4]{tom2008algebraic}.
\end{proof}

Note that the one-point space $\{*\}$ has a unique closure operator which is topological, and since our singular homology groups extend singular homology groups of topological spaces we have $H_k(*)=0$ for $k\ge 1$ and $H_0(*)=\mathbb{Z}$. This fact and \cref{theorem:homotopic_maps_induce_same_homomorphism} imply the following.

\begin{corollary}
Let $X$ be a contractible closure space. Then $H_k(X)=0$ for $k\ge 1$ and $H_0(X)=\mathbb{Z}$. \qed
\end{corollary}

We also associate homology groups to a pair of closure spaces $(X,A)$ via relative chains.
Let $X$ be a closure space and let $A$ be a subspace. Let $C_n(X,A)$, be the quotient group $C_n(X)/C_n(A)$. The boundary map $\partial_n:C_n(X)\to C_{n-1}(X)$ takes $C_n(A)$ to $C_{n-1}(A)$, hence it induces a quotient boundary map $\partial_n:C_n(X,A)\to C_{n-1}(X,A)$. Therefore we have a chain complex 
\[\cdots \to C_n(X,A)\xrightarrow{\partial}C_{n-1}(X,A)\to\cdots\]
The \emph{relative homology groups} $H_n(X,A)$, are the homology groups of the above chain complex. Standard arguments in homological algebra give the following.

\begin{theorem}{\cite[Theorem 5.15]{bubenik2021eilenbergsteenrod}}
\label{theorem:long_exact_sequence_of_relative_homology}
Given a closure space pair $(X,A)$, there exists a long exact sequence of homology groups:
\[\cdots\to H_n(A)\xrightarrow{i_*}\to H_n(X)\xrightarrow{j_*} H_{n}(X,A)\xrightarrow{\partial}H_{n-1}(A)\xrightarrow{i_*}H_{n-1}(X)\to\cdots \to H_0(X,A)\to 0\]
where $i:C_n(A)\to C_n(X)$ is the inclusion and $j:C_n(X)\to C_n(X,A)$ is the quotient map. \qed
\end{theorem}

\section{Hurewicz theorem in dimension $1$ for closure spaces}
\label{section:hurewicz}
Here we prove the Hurewicz theorem in dimension $1$ for the fundamental and singular homology groups of closure spaces we've been considering. 

Let $(X,x_0)$ be a pointed closure space. We define a function $h:\pi_1(X,x_0)\to H_1(X)$, called the \emph{Hurewicz map} in dimension 1. Let $[f]\in \pi_1(X,x_0)$ with a representative $f:(I,\partial I)\to (X,x_0)$. Since $I$ is homeomorphic to $|\Delta^1|$ we can think of $f$ as being a singular $1$-simplex in  $X$. Furthermore, since $f$ is a loop based at $x_0$, it follows that $f$ is also a cycle in $C_1(X)$ as $\partial(f)=f(1)-f(0)=x_0-x_0=0$. Set $h([f])$ to be the homology class of $f$.

\begin{wrapfigure}{5}{5.5cm}
\vspace{-10pt}
\centering
\begin{tikzpicture}[scale=1.1]
  \draw[thick] (0,0) node[below left] {$v_0$} -- (2,0) node[below right] {$v_1$} -- (2,2) node[above right] {$v_3$} -- (0,2) node[above left] {$v_2$} -- cycle;

  \draw[thick] (0,0) -- (2,2);

  \draw[-latex,thick] (1,0) -- (1.1,0);  
  \draw[-latex,thick] (2,1) -- (2,1.1);  
  \draw[-latex,thick] (1,2) -- (1.1,2);  
  \draw[-latex,thick] (0,1) -- (0,1.1);  
  \draw[-latex,thick] (0,0) -- (1.1,1.1);  

  \node at (1,-0.3) {$f_1$};
  \node at (1,2.3) {$f_2$};

  \node at (0.5,1.5) {$\sigma_1$};
  \node at (1.5,0.5) {$\sigma_2$};

\end{tikzpicture}
\caption{\small Homotopy $H$ as a sum of two singular $2$-simplices.}
\label{fig:h_is_well_defined}
\end{wrapfigure}
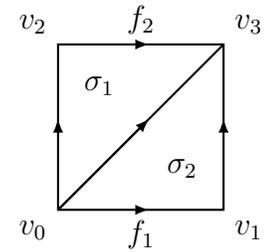
We show that $h$ is a well-defined function, i.e. $h([f])$ is always the same homology class that doesn't depend on the choice of representative $f$. 
Let $f_1, f_2$ be representatives of $[f]\in \pi_1(X,x_0)$. Then, there exists a homotopy (fixing endpoints) $H:I\times I\to X$ such that $H(t,0)=f_1(t)$ and $H(t,1)=f_2(t)$. Denote the vertices of $I\times I$ by $v_i$, $0\le i\le 3$ (\cref{fig:h_is_well_defined}). Then $H$ can be regarded as a sum of two singular $2$-simplices, $\sigma_1=[v_0,v_2,v_3]$ and $\sigma_2=[v_0,v_1,v_3]$. Let $c_{x_0}:I\to X$ be the constant singular $1$-simplex that represents the constant map at $x_0$. The boundary computations yield

\begin{gather*}
    \partial(\sigma_1)=H|_{[v_2,v_3]}-H|_{[v_0,v_3]}+H|_{[v_0,v_2]}=f_2-H|_{[v_0,v_3]}+H|_{[v_0,v_2]}\\
    \partial(\sigma_2)=H|_{[v_1,v_3]}-H|_{[v_0,v_3]}+H|_{[v_0,v_1]}=H|_{[v_1,v_3]}-H|_{[v_0,v_3]}+f_1
\end{gather*}

Therefore $\partial(\sigma_1-\sigma_2)=f_2-f_1+H|_{[v_0,v_2]}-H|_{[v_1,v_3]}=f_2-f_1+2c_{x_0}$. The constant singular simplex $c_{x_0}$ is a boundary of the constant map $|\Delta^2|\to X$ to $x_0$, hence it follows that $f_2-f_1$ is a boundary. Therefore, $f_1$ and $f_2$ are in the same homology class.

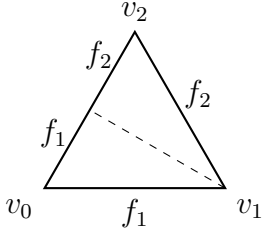
\begin{wrapfigure}[12]{l}{5.5cm}
\vspace{-10pt}
\centering
    \begin{tikzpicture}[scale=0.8]
  \def\sidelength{3}
  
  \pgfmathsetmacro\height{sqrt(3)/2*\sidelength}
  
  \draw[thick] (0,0) node[below left] {$v_0$} -- 
               (\sidelength,0) node[below right] {$v_1$} -- 
               (0.5*\sidelength,\height) node[above] {$v_2$} -- cycle;

  \draw[dashed] (\sidelength,0) -- (0.25*\sidelength,0.5*\height);
  
  \node at (0.05*\sidelength, 0.35*\height) {$f_1$};
  \node at (0.3*\sidelength, 0.85*\height) {$f_2$};
  \node at (0.85*\sidelength, 0.6*\height) {$f_2$};
  \node[below] at (0.5*\sidelength, 0) {$f_1$};
  
\end{tikzpicture}
\caption{\small The singular $1$-chain $f_2-f_1\cdot f_2+f_1$ as a boundary of a singular $2$-simplex.}
\label{fig:h_is_a_group_morphism}
\end{wrapfigure}
We now show that $h$ is more than a function, it is a group homomorphism $h:\pi_1(X,x_0)\to H_1(X)$. Let $[f_1],[f_2]\in \pi_1(X,x_0)$. We need to show that $h([f_1\cdot f_2])=h([f_1])+h([f_2])$. Let $\sigma:|\Delta^{2}|\to X$ be the orthogonal projection of $|\Delta^2|=[v_0,v_1,v_2]$ onto $[v_0,v_2]$ and is then followed by the concatenation $f_1\cdot f_2:I\to X$ (\cref{fig:h_is_a_group_morphism}). Therefore 
\[
\makebox[\textwidth][r]{$\partial(\sigma) = \sigma|_{[v_1,v_2]}-\sigma|_{[v_0,v_2]}+\sigma|_{[v_0,v_1]}=f_2-f_1\cdot f_2+f_1.$}\]
Thus, as $h([f_1])+h([f_2])$ and $f_1+f_2$ are in the same homology class, we have that $f_1\cdot f_2$ is in the same homology class as $h([f_1])+h([f_2])-\partial\sigma$. Hence, $h([f_1\cdot f_2])=h([f_1])+h([f_2])$.

\vspace{25pt}

Since $h$ is a group homomorphism, it follows that $h[\overline{f}]$ is in the same homology class as $-f$.

We now prove the Hurewicz theorem in dimension one for closure spaces. The proof uses the same ideas as the proof for topological spaces, such as the one presented in \cite[Theorem 2A.1]{MR1867354}. However, there is one step in the proof that does not extend directly to closure spaces, and thus the proof is not entirely trivial. It is~\cref{theorem:cw_complexes} that allows us to avoid the potential problem. Furthermore, even though we proved some properties of the map $h$ with respect to a basepoint, by \cref{lemma:fundamental_group_path_connected} if the space is path-connected the fundamental group does not depend on the basepoint. Therefore, in the statement of our theorem we consider $\pi_1(X)$.

\begin{theorem}[Hurewicz theorem in dimension $1$]
\label{theorem:homology_is_abelianization_of_fundamental_group}
 If $X$ is path-connected $h$ is surjective and its kernel is the commutator subgroup of $\pi_1(X)$. Hence $H_1(X)$ is the abelianization of $\pi_1(X)$. 
\end{theorem}

\begin{proof}
We first show that $h$ is surjective. Let $\sigma=\sum_{i}n_i\sigma_i$ be a $1$-cycle representing an element of $H_1(X)$. Splitting up multiples of a given $\sigma_i$ and relabeling we can assume that each $n_i=\pm 1$. Furthermore, we can always assume $n_i=1$ for all $i$ as we can reverse each path $\sigma_i$, and therefore we write $\sigma=\sum_{i}\sigma_i$ for the chosen $1$-cycle, without loss of generality.
If there exists an $i$ such that $\sigma_i$ is not a loop, then there exists a $j$ such that the concatenation $\sigma_i\cdot \sigma_j$ is defined, otherwise $\partial\sigma$ would not be $0$. Since $h$ is a group homomorphism, $\sigma_i+\sigma_j$ is in the same homology class as $\sigma_i\cdot\sigma_j$ and thus we can replace $\sigma_i+\sigma_j$ in the sum by the single term $\sigma_i\cdot\sigma_j$. By iterating this argument, we can assume withot loss of generality that each $\sigma_i$ is a loop.

Since $X$ is path-connected, for all $i$, there exists a path $\alpha_i$ from $x_0$ to the basepoint at $\sigma_i$. As $h$ is a group homomorphism, $\alpha_i\cdot \sigma_i\cdot \overline{\alpha}_i$ is in the same homology class as $\alpha_i+\sigma_i-\overline{\alpha}_i=\sigma_i$. Therefore, we can assume without loss of generality that each $\sigma_i$ is a loop at $x_0$. Combine all the $\sigma_i$ by concatenating them into a single loop based at $x_0$. Since $h$ is a group homomorphism, this single loop is in the same homology class as $\sigma$, showing that $h$ is surjective.

We now show that $\ker h$ is the commutator subgroup of $\pi_1(X)$. The commutator subgroup of $\pi_1(X)$ is contained in $\ker h$ as $H_1(X)$ is abelian. For the reverse inclusion we show that every homotopy class $[f]$ in the kernel of $h$ is trivial in the abelianization, $\pi_1(X)_{\text{ab}}$, of $\pi_1(X)$.

Let $[f]\in \ker h$, with a representative $f$. As a loop, $f$ is a cycle and since it is homologous to $0$, $f$ is the boundary of a singular $2$-chain $\sum_{i}n_i\sigma_i$. Without loss of generality, we assume each $n_i$ is $\pm 1$. For each $\sigma_i$, apply the boundary formula to write $\partial \sigma_i=\tau_{i0}-\tau_{i1}+\tau_{i2}$ for singular $1$-simplices $\tau_{ij}$. We thus have 
\[f=\partial (\sum_i n_i\sigma_i)=\sum_{i} n_i\partial \sigma_i=\sum_i n_i(\tau_{i0}-\tau_{i1}+\tau_{i2})=\sum_{i,j}(-1)^{j}n_i\tau_{ij}.\]
Since $f$ is a singular $1$-simplex, and $C_1(X)$ is the free abelian group generated by the singular $1$-simplices, the above equality together with $n_i=\pm 1$ implies $\tau_{ij}$ must form cancelling pairs, all except one which equals $f$. Glue together all the simplices that are the domains of $\sigma_i$'s, $|\Delta^2_i|$, by identifying cancelling pairs of edges and preserving orientations to get a $\Delta$-complex $K$. As $\Delta$-complexes are CW complexes, the closure operator for $K$ is topological by \cref{prop:limits_and_colimits_of_closure_spaces} and this is the same topological closure we would have if we had glued the $2$-simplices to get $K$ in $\cat{Top}$.

The maps $\sigma_i$ fit together to give a map $\sigma:K\to X$ (\cref{theorem:pasting_lemma_for_closure_spaces}). Deform $\sigma$, staying fixed on the edge corresponding to the domain of $f$ so that each vertex of $K$ maps to $x_0$. This can be done by considering paths from images of the vertices of $K$ to $x_0$ in $X$. Sliding along these paths defines a homotopy on the $0$-skeleton of $K$ and the edge corresponding to the domain of $f$. Apply the homotopy extension property in \cref{lemma:hep} to extend this homotopy to all of $K$, which is possible since $K$ and its $0$-skeleton are a pair of topological Hausdorff spaces. Restricting $\sigma$ to the simplices $\Delta^2_i$, we have a new chain $\sum_i n_i\sigma_i$ with boundary equal to $f$ and with all $\tau_{ij}$'s now being loops at $x_0$ after the deformation of $\sigma$. 

Since $\pi_1(X)_{\text{ab}}$ is abelian, we can write $[f]=\sum_{i,j}(-1)^jn_i[\tau_{ij}]$. Since we have $\partial\sigma_i=\tau_{i0}-\tau_{i1}+\tau_{i2}$ and $h$ is a homomorphism, we have $[\partial\sigma_i]=[\tau_{i0}]-[\tau_{i1}]+[\tau_{i2}]$ and therefore $[f]=\sum_in_i[\partial \sigma_i]$. Since each $\sigma_i$ gives a nullhomotopy of the composed loop $\tau_{i0}-\tau_{i1}+\tau_{i2}$, it follows that $[f]=0$ in $\pi_1(X)_{\text{ab}}$.
\end{proof}

\section{Excision and Mayer-Vietoris theorems for closure spaces}
\label{section:mayer_vietoris}
Here we show the excision and Mayer-Vietoris theorems for the singular homology groups we've been studying. Our proofs mostly rely on ideas in \cite{MR1867354,tom2008algebraic}. We also assume the reader is familiar with some minimal homological algebra techniques, in particular chain homotopy equivalences of chain complexes.

\subsection{Cone construction}
\label{subsection:cone_construction}
Let $X$ be a contractible closure space, $x_0\in X$, and let $\varepsilon_n:C_n(X)\to C_n(X)$ be a chain map defined by $\varepsilon_n=0$ for $n>0$ and $\varepsilon_0(\sum_in_i\sigma_i)=(\sum_in_i)c_{x_0}$ where $c_{x_0}:|\Delta^0|\to X$ is the constant map at $x_0$. Let $H:X\times I\to X$ be a homotopy between the identity map on $X$ and the constant map at $x_0$. Then $H$ induces a chain homotopy $h_n:C_{n-1}(X)\to C_{n}(X)$ between $\varepsilon$ and the identity $\id{C_{\bullet}(X)}$.

Let $q:|\Delta^{n-1}|\times I\to |\Delta^n|$ be the map 
$((s_0,\dots, s_{n-1}) ,t)\mapsto (t,(1-t)s_0,\dots, (1-t)s_{n-1})$.
Note that $|\Delta^n|$ can be realized as a cone of $|\Delta^{n-1}|$ and the map $q$ is the quotient map associated to the cone construction. For a singular $(n-1)$-simplex $\sigma:|\Delta^{n-1}|\to X$, there exists a unique singular $n$-simplex $h(\sigma)=h\sigma:|\Delta^n|\to X$ that satisfies $H(\sigma\times \id{I})=h(\sigma)\circ q$, as $q$ is a quotient map. The homomorphisms $h$ is obtained by linearly extending.

Furthermore, from the construction, we have that the faces of $h\sigma$ are $(h\sigma)_i=h(\sigma_{i-1})$ for $i>0$ and $(h\sigma)_0=\sigma$. Here we used simpler notation to denote by $\tau_i$ the restriction of a singular $n$-simplex $\tau$ to $[e_0,\dots, \hat{e}_i,\dots, e_n]$. Therefore we have for $n>1$, 
\[\partial_n (h\sigma)=(h\sigma)_0-\sum_{i=1}^n(-1)^{i-1}(h\sigma)_i=\sigma-\sum_{j=0}^{n-1}(-1)^jh(\sigma_{j-1})=\sigma-h(\partial_{n-1}\sigma),\]
and $\partial(h\sigma)=\sigma-c_{x_0}$ for a singular $0$-simplex $\sigma$. These formulas are to be expected; geometrically $h\sigma$ is a cone over $\sigma$ with apex $x_0$, and its boundary should be the base $\sigma$ together with all the cones with bases the boundaries of $\sigma$. From all of these equalities, it follows that $\partial h+h\partial=\id{C_{\bullet}(X)}-\varepsilon$ and therefore $h$ is a chain homotopy.

\subsection{Linear singular simplices}
\label{subsection:linear_simplices}

Let $Y\subset \mathbb{R}^N$ be a convex set and let $v_0,\dots ,v_n\in Y$. Denote the linear singular simplex $\sigma:|\Delta^n|\to Y$, $\sum_i t_ie_i\mapsto \sum_i t_i v_i$, by $\sigma=[v_0,\dots, v_n]$. In other words, $\sigma$ is the restriction to $|\Delta^n|$ of linear map $\mathbb{R}^{n+1}\to \mathbb{R}^N$ that sends $e_i$ to $v_i$, for all $i$. The linear singular simplices generate a subgroup $LC_n(Y)\subset C_n(Y)$. Boundaries of linear singular simplices are again linear singular simplices and so we have a subchain complex, $LC_{\bullet}(Y)\subset C_{\bullet}(Y)$. For every $v\in Y$, we have a straight line deformation retraction $Y\times I\to Y$ defined by $(x,t)\mapsto (1-t)x+tv$. Applying the cone construction described above to $[v_0,\dots, v_n]$ we obtain $[v,v_0,\dots, v_n]$. Denote the chain homotopy associated to the straight line deformation retraction by $v\cdot:LC_n(Y)\to LC_{n+1}(Y)$, $c\mapsto v\cdot c$.  As observed above, we have for a linear singular chain $c\in LC_n(Y)$, 
\[\partial v\cdot c=\begin{cases}
c-v\cdot \partial c, & n>0,\\
c-\varepsilon(c)v,& n=0,
\end{cases}\]
where $\varepsilon:LC_0(D)\to \mathbb{Z}$ is the mapping $\sum_i n_i\sigma_i\mapsto \sum_in_i$.

\begin{remark}
In \cref{subsection:cone_construction} we have argued that a contracting homotopy to a point will induce a contracting chain homotopy on general singular chains. However, here we are considering a straight line deformation retraction and that is why it indeed induces a contracting chain homotopy on linear singular chains as well.
\end{remark}

\subsection{Barycentric subdivisions of linear simplices} 

The barycenter of the linear singular simplex $\sigma=[v_0,\dots ,v_n]$ is the point $\textnormal{bar}(\sigma)=\frac{1}{n+1}\sum_{i=0}^n v_i$. Define a chain map $S_n:LC_n(Y)\to LC_n(Y)$ by mapping a linear singular simplex $\sigma:|\Delta^n|\to Y$ to
\[S_n(\sigma)=\begin{cases}
\sigma, & n=0,\\
\textnormal{bar}(\sigma)\cdot S_{n-1}(\partial \sigma), & n>0,
\end{cases}\]
and extending linearly. In other words, $S$ is the barycentric subdivision operator, that given a linear singular simplex $\sigma$, it first subdivides the faces of $\sigma$ inductively, then constructs cones with the bases the subdivisions and with apex the barycenter of the face.

We have that $S$ is indeed a chain map, $S_n:LC_n(Y)\to LC_n(Y)$. We proceed by induction and we have that for a singular linear simplex $\sigma$, $n\ge 2$,
\begin{align*}
\partial S_n\sigma&=\partial\textnormal{bar}(\sigma)\cdot S_{n-1}\partial\sigma= &\\
&=S_{n-1}\partial\sigma- \textnormal{bar}(\sigma)\cdot \partial (S_{n-1}\partial\sigma)= & \textnormal{ (since } \partial\textnormal{bar}(\sigma)\cdot=\id{LC_{\bullet}(D)}-\textnormal{bar}(\sigma)\cdot \partial)\\
&=S_{n-1}\partial\sigma- \textnormal{bar}(\sigma)\cdot (S_{n-2}\partial\partial \sigma)= & \textnormal{ (since } \partial S_{n-1}=S_{n-2}\partial \textnormal { by induction on } n)\\
&=S_{n-1}\partial \sigma & \textnormal{ (since } \partial^2=0).
\end{align*}

For $n=1$, let $\sigma=[v_0,v_1]$ be a linear singular $1$-simplex. Then $\partial\sigma=v_1-v_0$. We have 
\begin{gather*}
\partial S_1\sigma=\partial\textnormal{bar}(\sigma)\cdot S_0\partial\sigma=\partial\textnormal{bar}(\sigma)\cdot \partial\sigma=\textnormal{bar}\cdot (v_1-v_0)=\\
=\partial([\textnormal{bar}(\sigma),v_1]-[\textnormal{bar}(\sigma),v_0])=v_1-\textnormal{bar}(\sigma)-v_0+\textnormal{bar}(\sigma)=\\
=v_1-v_0=\partial\sigma=S_0\partial\sigma,
\end{gather*}
as $S_0$ is the identity.

The following facts about the barycentric subdivision can be found in algebraic topology textbooks.

\begin{lemma}{\cite[Lemma 9.4.2]{tom2008algebraic}}
    The diameter $d(v_0,\dots ,v_n)$ of the linear singular simplex $\sigma=[v_0,\dots ,v_n]$ with respect to the Euclidean norm is the maximum of the $||v_i-v_j||$, $i,j\in \{0,\dots,n\}$.\qed
\end{lemma}

\begin{lemma}{\cite[Lemma 9.4.3]{tom2008algebraic}}
\label{lemma:diameter_of_iterated_barycentric_subdivision}
    Let $v_0,\dots v_n\in \mathbb{R}^m$. Then $S_n[v_0,\dots ,v_n]$ is a linear combination of linear singular simplices with diameter at most $\frac{n}{n+1}d(v_0,\dots ,v_n)$. \qed
\end{lemma}

\subsection{Barycentric subdivision is chain homotopy equivalent to the identity}
We now construct a chain homotopy $T$ on $LC_{\bullet}(Y)$ between $S$ and $\id{LC_{\bullet}(Y)}$. We first set $T_0:LC_{0}(Y)\to LC_{1}(Y)$ to be the $0$ map, $T_0=0$. Then, we define $T_n:LC_n(Y)\to LC_{n+1}(Y)$ inductively by setting
\[T_n\sigma=\textnormal{bar}(\sigma)\cdot (\sigma -T_{n-1}\partial\sigma),\]
for a linear singular simplex $\sigma$ and extending linearly.

For $n=0$, since $\partial_0=0$ and $T_0=0$, the equation $\partial T+T\partial=\id{LC_{\bullet}(Y)}-S$ holds since $S_0=\id{LC_0(Y)}$.
We proceed by induction and we have that for a linear singular simplex $\sigma$, $n\ge 1$,
\begin{align*}
\partial T_n\sigma&=\partial\textnormal{bar}(\sigma)\cdot(\sigma- T_{n-1}\partial\sigma)= &\\
&=\sigma-T_{n-1}\partial\sigma- \textnormal{bar}(\sigma)\cdot \partial (\sigma-T_{n-1}\partial\sigma)= &  (\partial\textnormal{bar}(\sigma)\cdot=\id{LC_{\bullet}(D)}-\textnormal{bar}(\sigma)\cdot \partial)\\
&=\sigma-T_{n-1}\partial\sigma- \textnormal{bar}(\sigma)\cdot (\partial\sigma-\partial T_{n-1}\partial\sigma)=\\
&=\sigma-T_{n-1}\partial\sigma- \textnormal{bar}(\sigma)\cdot (S_{n-1}\partial\sigma + T_{n-2}\partial\partial \sigma)= & (\textnormal {by induction on } n)\\
&=\sigma-T_{n-1} \partial \sigma - S_n\sigma & ( \partial^2=0 \textnormal{ and } S_n\sigma=\textnormal{bar}(\sigma)\cdot S_{n-1}\partial\sigma).
\end{align*}
Therefore $\partial T+T\partial=\id{LC_{\bullet}(Y)}-S$ which means that $T$ is a chain homotopy between the identity and $S$ on $LC_{\bullet}(Y)$.

\subsection{Barycentric subdivision of general singular simplices}

Here we extend the definitions of barycentric subdivision $S$ to general singular simplices, and we argue that this extension is still chain homotopic to the identity.
Note that  $\id{|\Delta^n|}:|\Delta^n|\to |\Delta^n|$ is a linear singular simplex, and with notation introduced in \cref{subsection:linear_simplices} we can write $\id{|\Delta^n|}=[e_0,\dots, e_{n}]$. Given a singular $n$-simplex on a closure space $X$, $\sigma:|\Delta^n|\to X$, we have the chain map $\sigma_{\#}:LC_{\bullet}(|\Delta|^n)\subset C_{\bullet}(|\Delta^n|)\to C_{\bullet}(X)$, and in particular $\sigma=\sigma_{\#}[e_0,\dots,e_n]$. For $n\ge 0$, define $S_n:C_{n}(X)\to C_n(X)$ by 
\[S_n\sigma=\sigma_{\#}S_n[e_0,\dots, e_n],\]
and extending linearly. Geometrically, we are barycentrically subdividing the domain of $\sigma$, and the resulting $S_n\sigma$ is the linear combination (with appropriate signs) of the restrictions of $\sigma$ to the simplices in the subdivision. Furthermore, the morphisms $S_n$, $n\ge 0$, are a chain map since for a singular $n$-simplex $\sigma$ we have
\begin{align*}
\partial S_n\sigma&=\partial\sigma_{\#}S_n[e_0,\dots, e_n]=\sigma_{\#}\partial S_n[e_0,\dots, e_n]=\sigma_{\#}S_{n-1}\partial [e_0,\dots, e_n]=\\
&=\sigma_{\#}S_{n-1}\sum_{i}(-1)^i[e_0,\dots, \hat{e}_i,\dots, e_n]=\sum_i(-1)^i\sigma_{\#}S_{n-1}[e_0,\dots, \hat{e}_i,\dots e_n]=\\
&=\sum_i (-1)^iS_{n-1}\sigma|_{[e_0,\dots, \hat{e}_i,\dots, e_n]}=S_{n-1}\sum_i (-1)^i\sigma|_{[e_0,\dots, \hat{e}_i,\dots, e_n]}=S_{n-1}\partial \sigma.
\end{align*}
Now we define a chain homotopy between $S$ and identity on $C_{\bullet}(X)$ in a similar way, namely 
\[T_n\sigma=\sigma_{\#}T_n[e_0,\dots, e_n].\]
Using the same calculation as above where we replace $S$ by $T$, we get that $\sigma_{\#}T_{n-1}\partial [e_0,\dots, e_n]=T_{n-1}\partial \sigma$. It then follows that $T$ is indeed a chain homotopy between $S$ and $\id{C_{\bullet}(X)}$ since we have
\begin{align*}
\partial T_{n}\sigma&=\partial\sigma_{\#}T_n[e_0,\dots, e_n]=\sigma_{\#}\partial T_n[e_0,\dots, e_n]=\\
&=\sigma_{\#}([e_0,\dots, e_n]-S_n[e_0,\dots, e_n]-T_{n-1}\partial[e_0,\dots, e_n])=\\
&=\sigma-S_n\sigma-\sigma_{\#}T_{n-1}[e_0,\dots, e_n]=\sigma-S_n\sigma-T_{n-1}\partial\sigma
\end{align*}
This shows that the equality $\partial T+T\partial=\id{C_{\bullet}X}-S$ holds and thus $T$ is a chain homotopy.

\subsection{Small singular chains}
For a closure space $X$, let $\mathcal{U}=\{U_j\}_{j\in J}$ be an interior cover of $X$, and let $C_n^{\mathcal{U}}(X)$ be the subgroup of $C_n(X)$ consisting of chains $\sum_{i}n_i\sigma_i$ such that each $\sigma_i$ has their image contained in some element of $\mathcal{U}$. We call the chains in $C_{n}^{\mathcal{U}}(X)$ the $\mathcal{U}$-\emph{small singular chains} on $X$. The boundary map $\partial:C_n(X)\to C_{n-1}(X)$ takes $C_n^{\mathcal{U}}(X)$ to $C_{n-1}^{\mathcal{U}}(X)$, so the groups $C_n^{\mathcal{U}}(X)$ form a chain complex. The homology groups of this chain complex will be denoted by $H_n^{\mathcal{U}}(X)$.

\begin{lemma}
\label{lemma:lebesgue_number_for_small_simplices}
    Let $\sigma:|\Delta^n|\to X$ be a singular $n$-simplex. Then, there exists a $k\in \mathbb{N}$ such that each simplex in the chain $S_n^k\sigma$ has an image contained in an element of $\mathcal{U}$, where $S_n^k$ is the $k$-fold iteration of $S_n$.
\end{lemma}

\begin{proof}
    By \cref{lemma:finite_cover_whose_image_is_contained_in_covering_system}, there exists a finite open cover $\mathcal{V}$ of $|\Delta^n|$ such that for all $U\in \mathcal{U}$, there exists a $V\in \mathcal{V}$ such that $\sigma(V)\subset U$. Let $\delta >0$ be the Lebesgue number of this open cover, i.e., $\delta$ is such that every subset of $|\Delta^n|$ with diameter less than $\delta$ is contained in an element of the open cover $\mathcal{V}$. The simplices of $S^k_n$ are obtained by an application of $\sigma$ to $S^k_n[e_0,\dots, e_n]$. By \cref{lemma:diameter_of_iterated_barycentric_subdivision}, the diameter of these simplices is at most $(\frac{n}{n+1})^kd(e_0,\dots, e_n)$. This number will eventually be smaller than $\delta$ for $k$ large enough, and thus guaranteeing the statement of the lemma.
\end{proof}

We now prove that the chain groups $C_n(X)$ and $C^{\mathcal{U}}_n(X)$ are chain homotopy equivalent. Almost all of the arguments are borrowed from the proof of \cite[Proposition 2.21]{MR1867354} for the case of topological spaces and it is \cref{lemma:lebesgue_number_for_small_simplices} that allows us to extend them to all closure spaces.

\begin{proposition}
\label{prop:homology_of_covering_system}
The inclusion map $\iota:C_n^{\mathcal{U}}(X)\to C_n(X)$ is a chain homotopy equivalence. Hence $\iota$ induces isomorphisms $H_n^{\mathcal{U}}(X)\approx H_n(X)$ for all $n\in \mathbb{Z}$.
\end{proposition}

\begin{proof}
Consider the chain homotopy $T$ between $S$ and $\id{C_{\bullet}(X)}$, that was constructed in the previous subsection.
A chain homotopy between $\id{C_{\bullet}(X)}$ and the iterate $S^m$ is then given by $D_m=\sum_{0\le i < m} TS^i$, where we set $S^0=\id{C_{\bullet}(X)}$, as it follows that: 
\begin{gather*}
\partial D_m+D_m\partial=\sum_{0\le i < m} (\partial TS^i+TS^i\partial)=\sum_{0\le i < m}(\partial TS^i+T\partial S^i)=\\
\sum_{0\le i <m}(\partial T+T\partial)S^i=\sum_{0\le i <m}(S^i-S^{i+1})=S^0-S^m=\id{C_{\bullet}(X)}-S^m.
\end{gather*}
Define $D:C_n(X)\to C_{n+1}(X)$ by setting $D\sigma=D_{m(\sigma)}\sigma$, for a singular simplex $\sigma$ and extending linearly, where $m(\sigma)$ is the smallest $m$ such that $S^m\sigma$ is in $C_n^{\mathcal{U}}(X)$. Such an $m(\sigma)$ exists by \cref{lemma:lebesgue_number_for_small_simplices}.

Our goal is to find a chain map $\rho:C_n(X)\to C_n(X)$, with image in $C_n^{\mathcal{U}}(X)$, such that $D$ is a chain homotopy between $\rho$ and the identity. The image of $\rho$ being contained in $C_n^{\mathcal{U}}(X)$ and the chain homotopy will imply that the inclusion $\iota:C_n^{\mathcal{U}}(X)\to C_n(X)$ has a homotopy inverse, namely $\rho$. In other words, we want $\rho$ to be such that we have 
\[\partial D+D\partial=\id{C_{\bullet}(X)}-\rho.\]
We let the above equation define $\rho$ for us, that is we set $\rho=\id{C_{\bullet}(X)}-\partial D - D\partial$. We check that 
\begin{gather*}
	\partial\rho(\sigma)=\partial\sigma-\partial^2D\sigma-\partial D\partial \sigma=\partial\sigma-\partial D\partial\sigma,\\
	\rho(\partial\sigma)=\partial\sigma-\partial D\partial\sigma-D\partial^2\sigma=\partial\sigma-\partial D\partial \sigma.
\end{gather*}
which show that $\rho$ is indeed a chain map. Furthermore, $\rho$ takes chains in $C_n(X)$ to $\mathcal{U}$-small chains in $C_n^{\mathcal{U}}(X)$, as we have
\begin{align*}
\rho(\sigma)&=\sigma-\partial D\sigma-D\partial\sigma=\sigma-\partial D_{m(\sigma)}\sigma-D\partial\sigma=\\
&=S^{m(\sigma)}\sigma+D_{m(\sigma)}\partial\sigma-D\partial\sigma & \textnormal{ (since } \partial D_m+D_m\partial=\id{C_{\bullet}(X)}-S^m)
\end{align*}
By definition of $m(\sigma)$, $S^{m(\sigma)}\sigma$ is in $C_n^{\mathcal{U}}(X)$. The difference $D_{m(\sigma)}\partial\sigma-D\partial\sigma$ is a linear combination of terms of the form $D_{m(\sigma)}(\sigma_j)-D_{m(\sigma_j)}\sigma_j$ where $\sigma_j$ is a restriction of $\sigma$ to a face of $|\Delta^n|$. Note that $m(\sigma_j)\le m(\sigma)$ and thus $D_{m(\sigma)}(\sigma_j)-D_{m(\sigma_j)}\sigma_j$ consists of terms $TS^i(\sigma_j)$, $i\ge m(\sigma_j)$ and these terms also lie in $C^{\mathcal{U}}_n(X)$ as $T$ takes $C_{n-1}^{\mathcal{U}}(X)$ to $C_n^{\mathcal{U}}(X)$. Therefore, $\rho$ is a chain map $\rho:C_{\bullet}(X)\to C_{\bullet}^{\mathcal{U}}(X)$ and thus we have $\partial D+D\partial=\id{C_{\bullet}(X)}-\iota\rho$. Furthermore, $\rho\iota=\id{C_{\bullet}^{\mathcal{U}}(X)}$ as $D$ is $0$ on $C_{\bullet}^{\mathcal{U}}(X)$, as $m(\sigma)=0$ if $\sigma$ is $\mathcal{U}$-small and therefore the summation defining $D\sigma$ is empty. Therefore $\rho$ is a chain homotopy inverse for the inclusion $\iota:C_{\bullet}^{\mathcal{U}}(X)\to C_{\bullet}(X)$.
\end{proof}

Using \cref{prop:homology_of_covering_system}, we prove the following results.

\begin{theorem}[Excision]
\label{theorem:excision}
Let $X$ be a closure space. Suppose $\{A,B\}$ is an interior cover of $X$. Then there is an isomorphism for all $n\in \mathbb{Z}$, $H_n(X,A)\approx H_n(B,A\cap B)$. Equivalently, given subsets $Z\subset A\subset X$ such that the $c(Z)\subset i(A)$, the inclusion $(X-Z,A-Z)\to (X,A)$ induces isomorphisms $H_n(X-Z,A-Z)\approx H_n(X,A)$ for all $n$.
\end{theorem}

\begin{proof}
    To see that the two statements are equivalent, set $B=X\setminus Z$. Then it follows that $A\cap B=A\setminus Z$. Furthermore, by definition $i(B)=X\setminus c(X\setminus B)=X\setminus c(Z)$ which implies $c(Z)=X\setminus i(B)\subset i(A)$. Therefore $X=i(A)\cup i(B)$. We thus only prove the first statement about an interior $\{A,B\}$ cover of $X$.

    Let $\mathcal{U}$ denote the interior cover $\{A,B\}$ of $X$. At the end of the proof of \cref{prop:homology_of_covering_system} we showed the equality $\partial D+D\partial=\id{C_{\bullet}(X)}-\iota\rho$. Furthermore, all these maps take chains in $A$ to chains in $A$ and thus they induce quotient maps when factoring out chains in $A$. In particular, the formula $\partial D+D\partial=\id{C_{\bullet}(X)}-\iota\rho$ also holds for the inclusion $\iota:C_n^{\mathcal{U}}(X)/C_n(A)\to C_n^{\mathcal{U}}(X)/C_n(A)$. Therefore, $\iota:C_n^{\mathcal{U}}(X)/C_n(A)\to C_n^{\mathcal{U}}(X)/C_n(A)$ as a chain homotopy equivalence induces an isomorphism in homology, $H_n^{\mathcal{U}}(X,A)\approx H_n(X,A)$. Furthermore, the map $C_n(B)/C_n(A\cap B)\to C_n^{\mathcal{U}}(X)/C_n(A)$ induced by inclusion is an isomorphism since both quotients are free abelian groups with basis the singular simplices in $B$ that don't lie in $A$. Therefore, we have $H_n(X,A)\approx H_n(B,A\cap B)$.
\end{proof}

\begin{theorem}[Mayer-Vietoris]
\label{theorem:Mayer_Vietoris}
Let $X$ be a closure space and let $\{A,B\}$ be a interior cover of $X$. Consider the respective subspace closure operations on $A,B$ and $A\cap B$. Then there is a long exact sequence:
\[\cdots\to H_n(A\cap B)\to H_n(A)\oplus H_n(B)\to H_n(X)\to H_{n-1}(A\cap B)\to \cdots \]
\end{theorem}

\begin{proof}
    Let $\mathcal{U}$ denote the interior cover $\{A,B\}$. Then, by \cref{prop:homology_of_covering_system} there is an isomorphism $H_n(X)\approx H_n^{\mathcal{U}}(X)$. Furthermore, there is a short exact sequence of chain complexes
    \[0\to C_n(A\cap B)\xrightarrow{\varphi} C_n(A)\oplus C_n(B)\xrightarrow{\psi} C_n^{\mathcal{U}}(X)\to 0,\]
    given by $\varphi(x)=(x,-x)$ and $\psi(x,y)=x+y$. It follows from the definitions of $\varphi$ and $\psi$ that the above is indeed a short exact sequence of chain complexes. This short exact sequence induces the desired long exact sequence in homology, with the identification $H_n(X)\approx H_n^{\mathcal{U}}(X)$.
\end{proof}

Recall the notion of good pairs of closure spaces, $(X,A)$ (\cref{section:homotopy}). Using excision one can get the following familiar results applying verbatim arguments as in the case for topological spaces. 

\begin{proposition}
\label{prop:relative_and_reduced_absolute_homology}
For good pairs $(X,A)$, the quotient map $q:(X,A)\to (X/A,A/A)$ induces isomorphisms $q_*:H_n(X,A)\to H_n(X/A,A/A)\approx \tilde{H}_n(X/A)$ for all $n\in \mathbb{Z}$. 
\end{proposition}

\begin{proof}
The arguments in the proof of \cite[Proposition 2.22]{MR1867354} apply verbatim.  
\end{proof}

\begin{theorem}
\label{theorem:long_exact_sequence_reduced_homology}
Let $X$ be a closure space and let $A\subset X$ be non-empty and closed and suppose $(X,A)$ is a good pair. Then there is an exact sequence
\[\cdots \to \tilde{H}_n(A)\xrightarrow{i_*} \tilde{H}_n(X)\xrightarrow{q_*} \tilde{H}_n(X/A)\xrightarrow{\partial} \tilde{H}_{n-1}(A)\xrightarrow{i_*} \tilde{H}_{n-1}(X)\to\cdots\]
where $i$ is the inclusion $i:A\to X$, $q$ is the quotient map $q:X\to X/A$. 
\end{theorem}

\begin{proof}
    The result follows from \cref{theorem:long_exact_sequence_of_relative_homology,prop:relative_and_reduced_absolute_homology}.
\end{proof}

These results assemble together to give an axiomatic formulation of this singular homology for closure spaces in terms of Eilenberg-Steenrod type axioms~\cite{bubenik2021eilenbergsteenrod}. More specifically, the singular homology groups of closure spaces satisfy all the Eilenberg-Steenrod axioms for a homology theory. Whether these axioms guarantee uniqueness of homology theories in the category of closure spaces, or at least when restricted to some interesting full subcategory, remains open. Next we show that the suspension construction in $\cat{Cl}$ shifts homology groups by one degree. 

Given a closure space $(X,c)$, let $(\Sigma X,c_q)$ be the suspension of $(X,c)$, namely $\Sigma X=X\times I/\sim$ where $(x,0)\sim(y,0)$ and $(x,1)\sim (y,1)$  for all $x,y\in X$ and $X\times I$ has the product closure operation $c\times \tau$ while by $q$ we denote the quotient map and by $c_q$ the quotient closure operation.

\begin{proposition}
\label{prop:suspension}
Let $(X,c)$ be a closure space. Then $H_{k+1}(\Sigma X)\approx H_k(X)$ for all $k\ge 1$.
\end{proposition}

\begin{proof}
Let $\frac{1}{2}>\epsilon>0$ and let $A=\{[x,t]\in \Sigma X\,|\, t\in [0,\frac{1}{2}+\epsilon]\}$, $B=\{[x,t]\in \Sigma X\,|\, t\in [\frac{1}{2}-\epsilon,1]\}$. Observe that by definition of coequalizers in $\mathbf{Cl}$ we have 
\begin{align*}
i_{c_q}(A)&=\Sigma X-c_q(\Sigma X-A)=\Sigma X-q(c\times \tau(q^{-1}(\{[x,t]\in \Sigma X\,|\,t\in (\frac{1}{2}+\epsilon,1]\})))=\\
&=\Sigma X-q(c\times \tau (X\times (\frac{1}{2}+\epsilon,1]))=\Sigma X-q(c(X)\times \tau((\frac{1}{2}+\epsilon,1]))=\\
&=\Sigma X-q(X\times [\frac{1}{2}+\epsilon,1])
=\Sigma X-\{[x,t]\in \Sigma X\,|\, t\in [\frac{1}{2}+\epsilon,1]\}=\\
&=\{[x,t]\in \Sigma X\,|\,t\in [0,\frac{1}{2}+\epsilon)\}.
\end{align*}
Similarly, $i_{c_q}(B)=\{[x,t]\in \Sigma X\,|\, t\in (\dfrac{1}{2}-\epsilon,1]\}$ and thus $\{A,B\}$ is an interior cover as $X=\text{i}_{c_q}(A)\cup \text{i}_{c_q}(B)$ . By \cref{theorem:Mayer_Vietoris}, there exists a long exact sequence 
\[\cdots \to H_{k+1}(\Sigma X)\to H_k(A\cap B)\to H_k(A)\oplus H_k(B)\to H_k(\Sigma X)\to \cdots\]
Note that $A$ and $B$ are contractible with respect to their subspace closure operations by \cref{prop:cone} as each is homeomorphic to a cone of $(X,c)$.
Furthermore $(A\cap B,c_{A\cap B})$ deformation retracts to $(X\times \{\frac{1}{2}\},c_{X\times \{\frac{1}{2}\}})$ which is homeomorphic to $(X,c)$. Indeed, note that $(A\cap B,c_{A\cap B})$ is homeomorphic to $(X\times [\frac{1}{2}-\epsilon,\frac{1}{2}+\epsilon],c\times \tau)$. Consider the straight line deformation retraction $H$ of $[\frac{1}{2}-\epsilon,\frac{1}{2}+\epsilon]$ onto $\{\frac{1}{2}\}$. Then $\mathbf{1}_X\times H: (X\times [\frac{1}{2}-\epsilon,\frac{1}{2}+\epsilon])\to (X\times \{\frac{1}{2}\},c)$ is a deformation retraction. Thus the above sequence becomes
\[\cdots 0 \to H_{k+1}(\Sigma X)\to H_k(X)\to 0\to\cdots\]
for $k\ge 1$. Hence the result follows.
\end{proof}

Finally, we state an application of the excision theorem to a wedge of spaces.
The \emph{wedge sum} $\vee X_{\alpha}$ of a family $\{X_{\alpha}\}_{\alpha\in A}$ of closure spaces is a quotient space obtained by identifying chosen points $x_{\alpha}\in X_{\alpha}$ to a single point.

\begin{proposition}
\label{prop:homology_of_wedge}
Let $\vee X_{\alpha}$ be a wedge sum of closure spaces formed at basepoints $x_{\alpha}\in X_{\alpha}$ where the pairs $(X_{\alpha},x_{\alpha})$ are good. Then the inclusions $\iota_{\alpha}:X_{\alpha}\to \vee X_{\alpha}$ induce an isomorphism
$\oplus_{\alpha}\iota_{\alpha*}:\oplus\tilde{H}_n(X_{\alpha})\to \tilde{H}_n(\vee X_{\alpha})$ for all $n$.
\end{proposition}

\begin{proof}
The proof is a verbatim application of standard arguments as in the case of topological spaces, see for example \cite[Proposition 10.2.3]{tom2008algebraic}.
\end{proof}

\section{Computations}
\label{section:applications}
Here we apply the results of \cref{section:mayer_vietoris,section:hurewicz} and from \cite{rieser2021vcech} to compute examples of singular homology groups of closure spaces. We consider the circle $S^1$ with the geodesic metric, $d$. The closure operations we consider for $S^1$ will be induced by the geodesic metric. In particular, for a parameter $r\in \mathbb{R}$, we consider the closure space $(S^1,c_r)$, where 
\[c_r(A)=\{x\in S^1\,|\, \textnormal{dist}(x,A)=\inf_{y\in A}d(x,y)\le r\}.\]

It will be notationally convenient if the distance between points is a fraction of $1$ instead of a fraction of $2\pi$, thus whenever we write $S^1$ we mean a circle with circumference $1$. Roots of unity with varying closure operations will also be of interest to us, in particular
let $(\mathbb{Z}_n,c_{m})$ be the following closure space. We consider the set $\mathbb{Z}_n:=\{0,1,\dots, n-1\}$, for $n\in \mathbb{N}$. Define a metric on $\mathbb{Z}_n$, by $d(i,j)=\min\{|k|\,|\, k=i-j \textnormal { mod(n)}\}$. Alternatively, if we consider $\mathbb{Z}_n\subset S^1$ as a finite sample of points representing the $n$-th roots of unity, then $d$ is the restriction of the geodesic metric on $S^1$ to $\mathbb{Z}_n$ (rescaled by a factor of $2\pi$). Given $m\ge 0$, the closure $c_{m}$ defined by:
\[
c_m(i)=\{j\in \mathbb{Z}_n\,|\, d(i,j)\le k\} \textnormal{ and }
c_m(A)=\bigcup\limits_{i\in A}c_m(i)
\]
for all $A\subset \mathbb{Z}_n$.
For every closure $c_m$ we have an associated interior operation $i_m$ defined by
 \[i_{m}(A)=\mathbb{Z}_n\setminus c_{m}(\mathbb{Z}_n\setminus A)\,,\]
 for all $A\subset \mathbb{Z}_n$. For every subset $A\subset \mathbb{Z}_n$, we will denote the subspace closure operation on $A$ by $c_{A,m}$, or simply by $c_A$ when $m$ is clear from the context. 
 
A proof strategy for calculating homology groups of roots of unity will be to apply the Mayer-Vietoris theorem. For that we need to understand interior covers of roots of unity. All elements of interior covers we will consider will be homeomorphic to the space $(J_n,c_m)$ for some $n,m\in \mathbb{N}$ (\cref{def:discrete_interval}).

\begin{definition}
\label{def:discrete_interval}
Let $(\mathbb{Z},c_m)$ be the closure space defined by $c_m(i)=\{j\in \mathbb{Z}\,|\, |i-j|\le m\}$ for all $i\in \mathbb{Z}$ and $c_m(A)=\bigcup_{i\in A}c_m(i)$. Let $(J_n,c_m)$ denote $\{0,\dots , n\}$ with the subspace closure structure induced from $(\mathbb{Z},c_m)$. If $m=1$ we simply write $J_n$.
\end{definition}

The following observations will be key ingredients in the proofs of our major results.

\begin{lemma}
\label{lemma:discrete_intervals_are_contractible}
For all $n,m> 0$, $(J_n,c_m)$ is contractible.
\end{lemma}

\begin{proof}
$(\mathbb{Z},c_m)$ is contractible (\cite[Lemma 4.49]{rieser2021vcech}). Observe that $(J_n,c_m)$ is a retract of $(\mathbb{Z},c_m)$. Indeed, the inclusion map $i:(J_n,c_m)\to (\mathbb{Z},c_m)$ has a left inverse $r:(\mathbb{Z},c_m)\to (J_n,c_m)$ defined by 
\[r(i)=\begin{cases}
    i, & i\in \{0,\dots,n\}\\
    0, & i<0\\
    n, & n<i
\end{cases},\]
for all $i\in \mathbb{Z}$. Thus the claim follows, as retracts of contractible spaces are contractible by \cref{lemma:retracts_of_contractible_are_contractible}.
\end{proof}

\begin{lemma}
\label{lemma:discrete_intervals_are_homotopy_equivalent}
For all $n,m\ge 1$, $n<N$, $(J_N,c_m)$ deformation retracts to $(J_{n},c_m)$.
\end{lemma}

\begin{proof}
    It is sufficient to show that $(J_{n+1},c_m)$ deformation retracts to $(J_{n},c_m)$. The general case will follow by induction. Let $\iota:(J_n,c_m)\to (J_{n+1},c_m)$ be the inclusion map and $r:(J_{n+1},c_m)\to (J_n,c_m)$ the retraction mapping $n+1$ to $n$ and keeping everything else fixed. Then $r$ and $\iota$ are clearly continuous and $r\iota=\id{(J_n,c_m)}$. It remains to show that $\iota r\sim \id{(J_{n+1},c_m)}$.
    Define $H:(J_{n+1},c_m)\times J_1\to (J_{n+1},c_m)$ by setting $H(x,0)=\id{(J_{n+1},c_m)}(x)$ and $H(x,1)=\iota r(x)$. Then $H$ is continuous. Let $f:I\to J_1$ be defined by $f(0)=0$ and $f(x)=1$ for all $x>0$. As $J_1$ is the topological space with the indiscrete topology, $f$ is continuous. Precompose $H$ with $\id{(J_{n+1},c_m)}\times f$ to get the desired homotopy.
\end{proof}

\subsection{Initial calculations}

We calculate $H_1$ of $(S^1,c_r)$ and $(\mathbb{Z}_n,c_m)$, by building upon the work about fundamental groups by Rieser.

\begin{theorem}{\cite[Theorems 4.48 and 4.52]{rieser2021vcech}}
\label{theorem:fundamental_group_of_circle}
Consider the space $(S^1,c_r)$. Then
\[\pi_1(S^1,c_r)=\begin{cases}
\mathbb{Z}, & 0\le r<\frac{1}{3}\\
0, & \frac{1}{2}\le r
\end{cases}.\]
\end{theorem}

\begin{theorem}{\cite[Corollary 4.50 and Theorem 4.53]{rieser2021vcech}}
\label{theorem:fundamental_group_roots_of_unity}
Consider the space $(\mathbb{Z}_n,c_m)$. Then
\[\pi_1(\mathbb{Z}_n,c_m)=\begin{cases}
\mathbb{Z}, & 3\le 3m<n\\
0, & \floor{\frac{n}{2}} \le m \text{ or } n=3, 1\le m
\end{cases}.\]
\end{theorem}

Rieser also showed the following, for a wedge of circles. 

\begin{theorem}{\cite[Theorem 4.54]{rieser2021vcech}}
\label{theorem:fundamental_wedge_of_circles}
Let $X=(S^1,c_{r_1})\vee (S^1,c_{r_2})$ for some $r_1,r_2>0$. Then 
\[\pi_1(X)=\begin{cases}
F_2,& r_1,r_2<\dfrac{1}{3}\\
\mathbb{Z},& r_i<\dfrac{1}{3},r_j\ge \dfrac{1}{2}, i\neq j\\
0,& r_1,r_2\ge \dfrac{1}{2}
\end{cases},\]
where $F_2$ denotes the free group on two generators.
\end{theorem}

By \cref{theorem:homology_is_abelianization_of_fundamental_group}
we immediately have the following results.

\begin{corollary}
\label{theorem:homology_group_of_circle}
Consider the space $(S^1,c_r)$. Then
\[H_1(S^1,c_r)=\begin{cases}
\mathbb{Z}, & 0\le r<\frac{1}{3}\\
0, & \frac{1}{2}\le r
\end{cases}.\]
\end{corollary}

\begin{corollary}
\label{theorem:homology_of_roots_of_unity}
Consider the space $(\mathbb{Z}_n,c_m)$. Then
\[H_1(\mathbb{Z}_n,c_m)=\begin{cases}
\mathbb{Z}, & 3\le 3m<n\\
0, & \floor{\frac{n}{2}} \le m \text{ or } n=3, 1\le m
\end{cases}.\]
\end{corollary}

By \cref{theorem:homology_is_abelianization_of_fundamental_group,theorem:fundamental_wedge_of_circles} we have the following.
\begin{theorem}
\label{theorem:homology_wedge_of_circles}
Let $X$ be as in \cref{theorem:fundamental_wedge_of_circles}. Then 
\[H_1(X)=\begin{cases}
\mathbb{Z}\oplus \mathbb{Z},& r_1,r_2<\dfrac{1}{3}\\
\mathbb{Z},& r_i<\dfrac{1}{3},r_j\ge \dfrac{1}{2}, i\neq j\\
0,& r_1,r_2\ge \dfrac{1}{2}
\end{cases}.\]
\end{theorem}

An analogous result holds for roots of unity, but the arguments are a bit more involved. 

\begin{theorem}
\label{theorem:homology_wedge_of_roots_of_unity}
For $X=(\mathbb{Z}_{n_1},c_{m_1})\vee (\mathbb{Z}_{n_2},c_{m_2})$,
\[H_1(X)=\begin{cases}
\mathbb{Z}\oplus \mathbb{Z},& 3\le 3m_i<n_i, i=1,2 \\
\mathbb{Z},& 3m_i<n_i,\floor{\frac{n_j}{2}} \le m_j \text{ or } n_j=3, 1\le m_j, i\neq j\\
0,& \floor{\frac{n_i}{2}} \le m_i \text{ or } n_i=3, 1\le m_i, i=1,2
\end{cases}.\]
\end{theorem}

\begin{proof}
Without loss of generality, let $0$ be the basepoint of the wedge.
Suppose that $3\le 3m_i<n_i$ for $i=1,2$. In this case, $(\mathbb{Z}_{n_i},0)$ is a good pair in $(\mathbb{Z}_{n_i},c_{m_i})$ for $i=1,2$. Indeed, $c_{m_i}(0)$ is a neighborhood of $0$ in $(\mathbb{Z}_{n_i},c_{m_i})$ that is homeomorphic to $(J_{2m_i},c_{m_i})$ and thus it deformation retracts to $0$ by \cref{lemma:discrete_intervals_are_contractible}, for $i=1,2$. Thus, the first claim follows by \cref{prop:homology_of_wedge,theorem:homology_of_roots_of_unity}.

The other cases are similar, noting that if $\floor{\frac{n_i}{2}} \le m_i \text{ or } n_i=3, 1\le m_i$ for some $i=1,2$, then $(\mathbb{Z}_{n_i},0)$ is also a good pair in $(\mathbb{Z}_{n_i},c_{m_i})$; the neighborhood that deformation retracts to $0$ is $\mathbb{Z}_{n_i}$ as $(\mathbb{Z}_{n_i},c_{m_i})$ is the space with the indiscrete topology.
\end{proof}

\begin{remark}
\label{remark:interior_cover_restriction}
In order to do further computations we will be relying on \cref{theorem:Mayer_Vietoris,theorem:excision}. There are considerable restrictions for a cover of $(\mathbb{Z}_n,c_{m})$ to be an interior cover, see \cref{example:cover_that_is_not_interior_cover}.
\end{remark}

\begin{example}
\label{example:cover_that_is_not_interior_cover}
Consider the space $(\mathbb{Z}_4,c_{1})$. Let $A=\{0,1,3\}$ and $B=\{1,2,3\}$. Then $i_{1}(A)=\{0\}$ and $i_{1}(B)=\{2\}$. Thus $\{A,B\}$ is a cover of $(\mathbb{Z}_4,c_{1})$, but  not an interior cover. Thus we cannot apply \cref{theorem:Mayer_Vietoris,theorem:excision} on the cover $\{A,B\}$. However, the collection $A$, $B$, $C=\{2,3,0\}$ and $D=\{3,0,1\}$ is an interior cover.
\end{example}

\begin{theorem}
\label{theorem:higher_homology_of_roots_of_unity}
For all $n\ge 4$ and $k\in \mathbb{N}$, $H_k(\mathbb{Z}_n,c_{1})\approx H_k(\mathbb{Z}_4,c_{1})$. 
\end{theorem}

\begin{proof}
Suppose that $n\ge 4$. Let $A,B\subset \mathbb{Z}_n$ be the sets $A=\{1,\dots , n-3\}$ and $B=\{0,\dots , n-2\}$. Note that $i_{1}(B)=A$, thus $B$ is a neighborhood of $A$ in $(\mathbb{Z}_n,c_1)$. Observe also that $B$ deformation retracts to $A$ in $(\mathbb{Z}_n,c_1)$. Indeed, as subspaces $A$ and $B$ are homeomorphic to $J_{n-4}$ and $J_{n-2}$, respectively. Thus, by \cref{lemma:discrete_intervals_are_homotopy_equivalent} $B$ deformation retracts to $A$. Therefore, $(\mathbb{Z}_n,A)$ is a good pair in $(\mathbb{Z}_n,c_{1})$. By \cref{theorem:long_exact_sequence_reduced_homology} there exists a long exact sequence:
\[\cdots \to \tilde{H}_k(A)\to \tilde{H}_k(\mathbb{Z}_n,c_1)\to \tilde{H}_k(\mathbb{Z}_n/A)\to \tilde{H}_{k-1}(A)\to\cdots\]
where $A$ and $\mathbb{Z}_n/A$ are given the subspace and quotient closures, respectively. As $A\approx J_{n-4}$, $A$ is contractible by \cref{lemma:discrete_intervals_are_contractible}. Hence $\tilde{H}_k(A)=0$ for all $k$. Therefore, the above long exact sequence gives us isomorphisms $\tilde{H}_k(\mathbb{Z}_n,c_{1})\approx \tilde{H}_k(\mathbb{Z}_n/A)$. Observe that by construction $\mathbb{Z}_n/A$ is a four point set and in the quotient closure (\cref{ex:coequalizer}), $(\mathbb{Z}_n/A)$ is homeomorphic to $(\mathbb{Z}_4,c_{1})$. Thus $H_k(\mathbb{Z}_n,c_{1})\approx H_k(\mathbb{Z}_4,c_{1})$, for all $k\in \mathbb{N}$.
\end{proof}

\subsection{Higher homology groups}
We calculate some higher homology groups of $(\mathbb{Z}_n,c_m)$.

\begin{lemma}
\label{corollary:fundamental_and_homology_groups_of_3_roots}
$H_k(\mathbb{Z}_n,c_{1})=0$ and $\pi_k(\mathbb{Z}_n,c_{1})=0$ for all $k\ge 1$, $n\le 3$.
\end{lemma}

\begin{proof}
Note that the closure space $(\mathbb{Z}_n,c_{1})$ is a topological space for $n\le 3$; it is the $n\le 3$ point space with the indiscrete topology and thus contractible. Thus the result follows.
\end{proof}

\begin{theorem}
\label{theorem:higher_homology_for_n_greater_than_6m}
Suppose $n\ge 6m,m\ge 1$. Then $H_k(\mathbb{Z}_n,c_m)=0$ for $k\ge 2$.
\end{theorem}

\begin{proof}
Suppose $n=6m+l$ for some $l\ge 0$. Consider the sets 
\[A=\{0,1,\dots, 5m+l\} \textnormal{ and } B=\{3m,3m+1,\dots, 6m+l-1,0,1,\dots, 2m-1\}.\]
Direct computation shows that 
\[i_{m}(A)=\{m,m+1,\dots, 4m+l\} \textnormal{ and } i_{m}(B)=\{4m,4m+1,\dots, 6m+l-1,0,1,\dots, m-1\}.\]
Thus as $\mathbb{Z}_n=i_{m}(A)\cup i_{m}(B)$, $\{A,B\}$ is an interior cover of $(\mathbb{Z}_n,c_m)$. By \cref{theorem:Mayer_Vietoris} there exists a long exact sequence:
\[\cdots \to H_k(A\cap B)\to H_k(A)\oplus H_k(B)\to H_k(\mathbb{Z}_n,c_m)\to \cdots, \]
where $A,B$ and $A\cap B$ have their respective subspace closures induced from $(\mathbb{Z}_n,c_m)$.
By construction $A\cap B=V_1\cup V_2$ where 
\[V_1=\{0,1,\dots, 2m-1\} \textnormal{ and } V_2=\{3m,3m+1,\dots, 5m+l\}.\]
Note that $c_{A\cap B}(V_1)=V_1$ and $c_{A\cap B}(V_2)=V_2$, and thus $V_1$ and $V_2$ are closed in $A\cap B$. Furthermore, they are each other's complement in $A\cap B$. Hence they form a separation of $(A\cap B,c_{A\cap B})$. Thus, as $I$ is connected, any continuous function $\gamma:I\to (A\cap B,c_{A\cap B})$ must lie entirely in $V_1$ or $V_2$. Hence, $(A\cap B,c_{A\cap B})$ has two path components, $V_1$ and $V_2$. Now observe that both $V_1$ and $V_2$ are contractible in $(A\cap B,c_{A\cap B})$. Indeed, $(V_1,c_{V_1})$ and $(V_2,c_2)$ are homeomorphic to $(J_{2m-1},c_m)$ and $(J_{2m+l},c_m)$, respectively, which are both contractible by \cref{lemma:discrete_intervals_are_contractible}.  In particular, $(A\cap B,c_{A\cap B})$ is contractible to the two point space with the discrete topology. Therefore $H_k(A\cap B)=0$ for $k\ge 1$. Thus, for $k\ge 2$, the long exact sequence from above yields: $H_{k+1}(\mathbb{Z}_n,c_{m})\approx H_k(A)\oplus H_k(B)$.
Note that $A$ and $B$ are homeomorphic to $(J_{5m+l},c_m)$ and $(J_{5m+l-1},c_m)$ and are thus contractible by \cref{lemma:discrete_intervals_are_contractible}.
Hence $H_k(\mathbb{Z}_n,c_{m})=0$, $k\ge 2$.
\end{proof}

\begin{theorem}
\label{theorem:higher_homology_for_n_greater_than_4m}
Suppose $n\ge 4m$. Then $H_k(\mathbb{Z}_n,c_{m})=0$ for $k\ge 2$.
\end{theorem}

\begin{proof}
If $m=0$, the space $(\mathbb{Z}_n,c_{m})$ is a topological space; it is the set of $n$ many points with the discrete topology. Thus $H_k(\mathbb{Z}_n,c_{m})=0$ for $k\ge 1$ and the result follows. 

Suppose $m\ge 1$ and $n=4m+l$ for some $l\ge 0$. 
Consider the sets
\[A=\{0,1,\dots, m+l-1\} \textnormal{ and  } B=\{3m+l,3m+l+1,\dots, 0,1,\dots, 2m+l-1\}.\]
Observe that $B$ is a neighborhood of $A$ in $(\mathbb{Z}_n,c_{m})$. Furthermore $B$ deformation retracts to $A$. To see this, realize that $B$ and $A$ as subspaces of $(\mathbb{Z}_n,c_m)$ are homeomorphic to $(J_{3m+l-1},c_m)$ and $(J_{m+l-1},c_m)$, respectively, and apply \cref{lemma:discrete_intervals_are_homotopy_equivalent}. 
Therefore $(\mathbb{Z}_n,A)$ is a good pair in $(\mathbb{Z}_n,c_m)$.  Thus, there exists a long exact sequence:
\[\cdots\to H_k(A)\to H_k(\mathbb{Z}_n,c_{m})\to \tilde{H}_k(\mathbb{Z}_n/A)\to\cdots,\]
where $A$ and $\mathbb{Z}/A$ are given the subspace and quotient closures, respectively. As $A$ is homemomorphic to $(J_{m+l-1},c_m)$, it is contractible by \cref{lemma:discrete_intervals_are_contractible}. Thus the long exact sequence above yields:
 \[H_k(\mathbb{Z}_n,c_{m})\approx \tilde{H}_k(\mathbb{Z}_n/A),\, k\ge 2.\]
Let $N\in \mathbb{N}$, $N\ge 6m$ and consider the closure space $(\mathbb{Z}_N,c_{m})$ and the sets 
\[V=\{0,1,\dots, N-3m-1\} \textnormal{ and } U=\{N-m,N-m+1,\dots, 0,1,\dots ,N-2m-1\}.\]
By construction $U$ is a neighborhood of $V$ in $(\mathbb{Z}_N,c_{m})$. Using similar arguments as when we showed $B$ deformation retracts to $A$, we can show $U$ deformation retracts to $V$. Then $(\mathbb{Z}_N,V)$ is a good pair in $(\mathbb{Z}_N,c_m)$.   Furthermore, $V$ is homeomorphic to $(J_{N-3m-1},c_m)$ and thus contractible by \cref{lemma:discrete_intervals_are_contractible}. Thus, from the long exact sequence for homology for the pair $(\mathbb{Z}_N,V)$ in $(\mathbb{Z}_N,c_m)$ for $k\ge 2$ we have isomorphisms
\[H_k(\mathbb{Z}_N,c_{m})\approx \tilde{H}_k(\mathbb{Z}_N/V),\]
where $\mathbb{Z}_N/V$ has the quotient closure. However, note that $\mathbb{Z}_N/V$ is homeomorphic to $\mathbb{Z}_n/A$, with respect to their quotient closure operators induced from $(\mathbb{Z}_N,c_m)$ and $(\mathbb{Z}_n,c_m)$, respectively. Thus for $k\ge 2$ we have 
\[H_k(\mathbb{Z}_N,c_{m})\approx \tilde{H}_k(\mathbb{Z}_N/V)\approx \tilde{H}_k(\mathbb{Z}_n/A)\approx H_k(\mathbb{Z}_n,c_{m}).\] 
By \cref{theorem:higher_homology_for_n_greater_than_6m} and our choice of $N$, we have that $H_k(\mathbb{Z}_n,c_{m})\approx 0$ for $k\ge 2$.
\end{proof}

By \cref{prop:suspension,theorem:homology_of_roots_of_unity,theorem:higher_homology_for_n_greater_than_4m}, we have the following.

\begin{proposition}
\label{prop:suspension_1}
If $n\ge 4m$, then $H_k(\sum (\mathbb{Z}_n,c_m))\approx\mathbb{Z}$ for $k=0,2$ and $0$ otherwise.
\end{proposition}

\subsection{Varying the unity structure}
\begin{figure}[H]
\centering
\begin{tikzpicture}
[
dot/.style={draw,fill,circle,inner sep=1pt},
squarednode/.style={rectangle, draw, fill, inner sep=1pt},
starnode/.style={star, draw, fill, inner sep=1pt},
diamondnode/.style={diamond, draw, fill, inner sep=1pt}
]
\foreach \i in {0,...,10} {

   \node[dot,label={\i*360/11-(\i==11)*45:$\i$}] (w\i) at ( \i*360/11:1)   {};
   }
   \node[dot,color=red,scale=2.5] at (0:1){};
   \node[dot,color=red,scale=2.5] at (360/11:1){};
   \node[squarednode,color=blue,scale=2.5] at (3*360/11:1){};
   \node[squarednode,color=blue,scale=2.5] at (4*360/11:1){};
   \node[starnode,color=green,scale=2.5] at (5*360/11:1){};
   \node[starnode,color=green,scale=2.5] at (6*360/11:1){};
   \node[diamondnode,color=pink,scale=2.5] at (8*360/11:1){};
   \node[diamondnode,color=pink,scale=2.5] at (9*360/11:1){};
\end{tikzpicture}
\quad\quad\quad\quad
\begin{tikzpicture}
[
dot/.style={draw,fill,circle,inner sep=1pt},
squarednode/.style={rectangle, draw, fill, inner sep=1pt},
starnode/.style={star, draw, fill, inner sep=1pt},
diamondnode/.style={diamond, draw, fill, inner sep=1pt}
]
\foreach \i in {0,...,6} {

   \node[dot,label={\i*360/7-(\i==7)*45:$\i$}] (w\i) at ( \i*360/7:1)   {};
   }
   \node[dot,color=red,scale=2.5] at (0:1){};
   \node[squarednode,color=blue,scale=2.5] at (2*360/7:1){};
   \node[starnode,color=green,scale=2.5] at (3*360/7:1){};
   \node[diamondnode,color=pink,scale=2.5] at (5*360/7:1){};
\end{tikzpicture}
\caption{The closure space $(\mathbb{Z}_{11},c_3)$ (left) and its quotient that is homeomorphic to $(\mathbb{Z}_{7},c_2)$ (right).}
\label{fig:quotients_preserving_homology}
\end{figure}
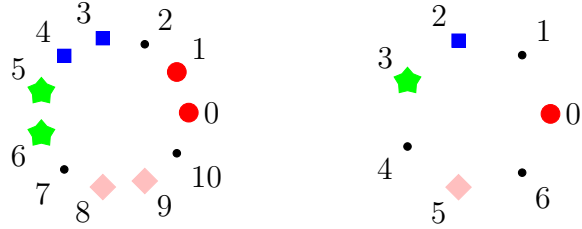

We investigate how different roots of unity structures are related in terms of homology groups. Consider \cref{fig:quotients_preserving_homology}, where we see that if we quotient $(\mathbb{Z}_{11},c_3)$ by identifying four suitable pairs of points, we recover the closure structure of  $(\mathbb{Z}_{7},c_2)$. If each of the four sets of interest give us a good pair in $(\mathbb{Z}_{11},c_3)$ the homology groups will be preserved and thus $H_k(\mathbb{Z}_{11},c_3)\approx H_k(\mathbb{Z}_7,c_2)$.

We generalize the idea in \cref{fig:quotients_preserving_homology} with the following results.

\begin{theorem}
\label{theorem:higher_homology_for_n_between_3m_and_4m_part_1}
For $m\ge 2$, suppose that $n=3m+l$, $\max(0, m-3)<l<m$. Then \[H_k(\mathbb{Z}_n,c_{m})\approx H_k(\mathbb{Z}_{n+4},c_{m+1}),\textnormal{ for all } k\in \mathbb{N}.\]
\end{theorem}

\begin{proof}
Let $A_1,A_2,A_3,A_4\subset \mathbb{Z}_{n+4}$ be the subspaces $A_1=\{0,1\}$, $A_2=\{m+1,m+2\}$, $A_3=\{2m+1,2m+2\}$ and $A_4=\{3m+2,3m+3\}$. Each $A_i$ has a neighborhood $B_i$ in $(\mathbb{Z}_{n+4},c_{m+1})$ that deformation retracts onto $A_i$, for $1\le i\le 4$. Indeed, if $B_i=c_{m+1}(A_i)$, then $i_{m+1}(B_i)=A_i$,  for $1\le i\le 4$. Furthermore, $(B_i,c_{B_i})$ is homeomorphic to $(J_{2m+3},c_{m+1})$ and thus deformation retracts to $(A_i,c_{A_i})\approx J_1$ by \cref{lemma:discrete_intervals_are_homotopy_equivalent}, for $1\le i\le 4$. Hence $(\mathbb{Z}_{n+4},A_1)$ is a good pair in $(\mathbb{Z}_{n+4},c_{m+1})$. If $q_1:\mathbb{Z}_{n+4}\to \mathbb{Z}_{n+4}/A_1$ is the quotient map,  by \cref{theorem:long_exact_sequence_reduced_homology} there is a long exact sequence
\[\cdots\to\tilde{H}_k(A_1)\to \tilde{H}_k(\mathbb{Z}_{n+4},c_{m+1})\to \tilde{H}_k(\mathbb{Z}_{n+4}/A_1,c_{q_1})\to \tilde{H}_{k-1}(A_1)\to \cdots\]
Since $(A_1,c_{A_1})\approx J_1$ is contractible by \cref{lemma:discrete_intervals_are_contractible}, the above sequence yields isomorphisms \[\tilde{H}_k(\mathbb{Z}_{n+4},c_{m+1})\approx \tilde{H}_k(\mathbb{Z}_{n+4}/A_1,c_{q_1}),\, \forall k\in \mathbb{N}.\]
As $q_1(B_2)$ is a neighborhood of $A_2$ that deformation retracts onto $q_1(A_2)$ in $(\mathbb{Z}_{n+4}/A_1,c_{q_1})$,
$(\mathbb{Z}_{n+4}/A_1,q_1(A_2))$ is a good pair in $(\mathbb{Z}_{n+4}/A_1,c_{q_1})$. By \cref{theorem:long_exact_sequence_reduced_homology} we get a long exact sequence
\[\cdots\to\tilde{H}_k(q_1(A_2))\to \tilde{H}_k(\mathbb{Z}_{n+4}/A_1)\to \tilde{H}_k((\mathbb{Z}_{n+4}/A_1)/q_1(A_2))\to \tilde{H}_{k-1}(q_1(A_2))\to \cdots\]
where we avoided specifying the closure operations of each space for simplicity, but they should be clear from context.
As $q_1(A_2)$ is contractible in $(\mathbb{Z}_{n+4}/A_1)/q_1(A_2)$, the long exact sequence yields isomorphisms 
\[\tilde{H}_k(\mathbb{Z}_{n+4}/A_1)\approx \tilde{H}_k((\mathbb{Z}_{n+4}/A_1)/q_1(A_2)),\, \forall k\in \mathbb{N}.\] 
Thus we have isomorphisms $\tilde{H}_k(\mathbb{Z}_{n+4},c_{m+1})\approx \tilde{H}_k((\mathbb{Z}_{n+4}/A_1)/q_1(A_2))$. Repeat the above arguments with sets $A_3$ and $A_4$, to get isomorphisms 
\[\tilde{H}_k(\mathbb{Z}_{n+4},c_{m+1})\approx \tilde{H}_k(\mathbb{Z}_{n+4}/\sim,c_{q}),\forall k\in \mathbb{N},\] 
where $\sim$ is the equivalence relation $k\sim j$ if $k,j\in A_i$, $i=\{1,2,3,4\}$, and $q$ is the corresponding quotient map. Note that there is a homeomorphism $f:(\mathbb{Z}_{n+4}/\sim,c_q)\to (\mathbb{Z}_n,c_m)$ given by 
\[f([i])=\begin{cases}
0, & i\in A_1\\
i-1, & i\in \{2,3,\dots, m\}\\
m, & i\in A_2\\
i-2, & i\in \{m+3,m+4,\dots ,2m\}\\
2m-1, & i\in A_3\\
i-3, & i\in \{2m+3, 2m+4,\dots ,3m+1\}\\
3m-1, & i\in A_4\\
i-4, & \ i\in \{3m+4,3m+5,\dots ,3m+l+3\}
\end{cases},\]
where $[i]$ is the equivalence class of $i$ in $\mathbb{Z}_{n+4}/\sim$.
As the sets $\mathbb{Z}_{n+4}/\sim$ and $\mathbb{Z}_n$ are finite, it suffices to show $f(c_q([i]))=c_{m}(f([i]))$ for all $[i]\in \mathbb{Z}_{n+4}/\sim$ for $f$ to be a homeomorphism. By \cref{ex:coequalizer} we have the following calculations:
\begin{itemize}[left=0pt]
\item Suppose $i\in A_1$. Then $c_q([i])=\{[2m+l+3],[2m+l+4],\dots [0],[1],\dots,[m+1],[m+2]\}.$
As $\max(0,m-3)< l<m$, we have $4\le l+3\le m+2$. Thus $2m+l+3\in\{2m+4,2m+3,\dots ,3m+2\}$. Therefore
\[f(c_q([i]))=\{2m+l,2m+l+1,\dots ,3m+l-1,0,1,\dots m-1,m\}=c_{m}(0)=c_{m}(f([i])).\]

\item Suppose $i\in \{2,3,\dots ,m\}$. Then $i+m+1\in \{m+3,m+4,\dots,2m+1\}$. Note that the distance between $2m+l+i+3$ and $i$ is $m+1$ in $\mathbb{Z}_{n+4}$. Since $\max(0,m-3)<l<m$, if $m=2$, then $0<l<m$. If $l=1$ and $i=2$ we have $2m+l+i+3=3m+4$. If $m\ge 3$, then $m-3<l<m$, and we have $2m+(m-3)+i+3=3m+i<2m+l+i+3$ and as $i\in \{2,3,\dots ,m+1\}$, $3m+i$ is at least $3m+2$. Hence \[c_q([i])=\{[2m+l+i+3],[2m+l+i+4],\dots ,[0],[1],\dots [i],[i+1],\dots ,[i+m],[i+m+1]\},\] 
and therefore
\[f(c_q([i]))=\{2m+l+i-1,2m+l+i,\dots, i-1,i,\dots, i+m-1 \}=c_{m}(i-1)=c_{m}(f([i])).\]

\item Suppose $i\in A_2$. Then $c_q([i])=\{[0],[1],[2],\dots ,[m+1],[m+2],\dots ,[2m+2],[2m+3]\}$,
and thus \[f(c_q([i]))=\{0,1,\dots, m-1,m,\dots ,2m-1,2m\}=c_m(m)=c_{m}(f([i])).\]

\item Suppose $i\in \{m+3,m+4,\dots,2m\}$. We have
\[c_q([i])=\{[i-m-1],[i-m],\dots ,[i],[i+1],\dots, [i+m], [i+m+1]\}.\]
As $i\in \{m+3,m+4,\dots ,2m\}$, we have $i-m-1\in \{2,3,\dots,m-1\}$ and $i+m+1\in \{2m+4,2m+5,\dots,3m+1\}$. Hence 
\[
    f([i-m-1])=i-m-1-1=i-m-2 \textnormal{ and }
    f([i+m+1])=i+m+1-3=i+m-2.
\]
Therefore \[
    f(c_q([i]))=\{i-m-2,i-m-1,\dots ,i-2,\dots, i+m-3,i+m-2 \}=c_{m}(i-2)=c_{m}(f([i])).
\]

\item Suppose $i\in A_3$. Then 
    $c_q([i])=\{[m],[m+1],\dots ,[2m+1],[2m+2],\dots ,[3m+2],[3m+3]\}$.
    Therefore
    \[
        f(c_q([i]))=\{m-1,m,\dots, 2m-1,2m, \dots, 3m-2,3m-1\}
        =c_{m}(2m-1)=c_{m}(f([i])).
    \]

\item Suppose $i\in\{2m+3,2m+4,\dots ,3m+1\}$. Note that the distance between $3m+1$ and $m-l-2$ is $m+1$ in $\mathbb{Z}_{n+4}$. Therefore $i+m+1\in \{3m+4,3m+5,\dots ,m-l-2\}$ (since, $\max(0,m-3)<l$, we have $m-l-2<1$).  Furthermore, we also have $i-m-1\in \{m+2,m+3,\dots,2m\}$. Since
\[c_q([i])=\{[i-m-1],[i-m],\dots [i],[i+1],\dots ,[i+m],[i+m+1]\},\]
we have 
\begin{align*}
    f(c_q([i]))=&\{i-m-3,i-m-2,\dots,i-3,i-2,\dots , i+m-4,i+m-3\}=\\
    =&c_{m}(i-3)=c_{m}(f([i])).
\end{align*}

\item Suppose $i\in A_4$. 
Note that the distance between $3m+3$ and $m-l$ is $m+1$ in $\mathbb{Z}_{n+4}$. Since $l<m$ implies $0< m-l<m$ we have
 \[c_q([i])=\{[2m+1],[2m+2],\dots ,[3m+1],[3m+2],\dots ,[0],\dots ,[m-l]\}.\] 
Therefore
\[
f(c_q([i]))=\{2m-1,2m,\dots ,3m-2,3m-1,\dots ,0,\dots ,m-l-1\}
=c_{m}(3m-1)=c_{m}(f([i])).    
\]

\item Suppose $i\in \{3m+4,3m+5,\dots,3m+l+3\}$. 
Then, we have $i-m-1\in \{2m+3,2m+4,\dots,2m+l+2\}$. Note that $l<m$ implies $2m+l+2<3m+2$, hence $i-m-1\in \{2m+3,2m+4,\dots,3m+1\}$. Thus, we have $f([i-m-1])=i-m-4$. Furthermore, the distance between $3m+4$ and $m-l+1$ is $m+1$ in $\mathbb{Z}_{n+4}$.  Observe that $l<m$ implies $m-l+1\ge 2$. Additionally, the distance between $3m+l+3$ and $m$ is $m+1$ in $\mathbb{Z}_{n+4}$.
We also have that 
\[i+m+1=i+m+1-3m-l-4 \textnormal{ mod(n+4)}=i-2m-l-3 \textnormal{ mod(n+4)}.\] Since $i-2m-l-3$ is between $2$ and $m$ as argued above, we have $f([i+m+1])=i-2m-l-3-1=i-2m-l-4$ by the definition of $f$. 
Therefore
\begin{align*}
    f(c_q([i]))=&\{i-m-4,i-m-3,\dots,i-4,i-3,\dots, 0,1\dots, i+m-1,i-2m-l-4\}=\\
    =&c_{m}(i-4)=c_{m}(f([i])).
\end{align*}
\end{itemize}
Thus $(\mathbb{Z}_{n+4}/\sim, c_q)$ is homeomorphic to $(\mathbb{Z}_n,c_{m})$ and therefore the result follows.
\end{proof}

\begin{theorem}
\label{theorem:higher_homology_for_n_between_3m_and_4m_part_2}
There exist isomorphisms
\[H_k(\mathbb{Z}_{4m-4},c_{m})\approx H_k(\mathbb{Z}_{4m},c_{m+1}),\, \textnormal{for } m\ge 5, k\in \mathbb{N}.\]
\end{theorem}

\begin{proof}
Let $m\ge 5$.
Let $A_1,A_2,A_3,A_4\subset \mathbb{Z}_{n+4}$ be the subsets $A_1=\{0,1\},A_2=\{m,m+1\}$, $A_3=\{2m,2m+1\}$ and $A_4=\{3m,3m+1\}$. We can argue like in the proof of \cref{theorem:higher_homology_for_n_between_3m_and_4m_part_1} that $(\mathbb{Z}_{4m},A_i)$ is a good pair in $(\mathbb{Z}_{4m},c_{m+1})$, for $1\le i\le 4$. Since we want to construct the sets $A_i$, $1\le i\le 4$, so that they don't overlap with one another, we need $\mathbb{Z}_{4m}$ to have at least $8$ points and thus we need $m\ge 2$. Furthermore, the argument that each $(\mathbb{Z}_{4m},A_i)$, $1\le i\le 4$ is a good pair is actually not true when $m=3,4$. For example if $m=3$, then $c_{m+1}(A_1)=c_4(A_1)=\{8,9,10,11,0,1,2,3,4,5\}$ but $5$ and $8$ are a distance less than $4$ from one-another and are hence in each other closures. Thus the subspace $c_4(A_1)$ is not homeomorphic to $(J_9,c_4)$ and so we can't argue it deformation retracts to $A_1$. Similarly for the case $m=4$. However, these arguments start working when $m\ge 5$. Furthermore, after taking quotients and collapsing each $A_i$ to a point we get isomorphisms $H_k(\mathbb{Z}_{4m},c_{m+1})\approx H_k(\mathbb{Z}_{4m}/\sim,c_q)$, where $c_q$ denotes the quotient closure. 
There is a homeomorphism $f:(\mathbb{Z}_{4m}/\sim,c_q)\to (\mathbb{Z}_{4m-4},c_{m})$, defined by

\[f([i])=\begin{cases}
0,& i\in A_1\\
i-1, & i\in \{2,3,\dots ,m-1\}\\
m-1, & i\in A_2\\
i-2, & i\in \{m+2,m+3\dots ,2m-1\}\\
2m-2, & i\in A_3\\
i-3, & i\in \{2m+2,2m+3,\dots ,3m-1\}\\
3m-3, & i\in A_4\\
i-4, & i\in \{3m+2,3m+3,\dots, 4m-1\}
\end{cases},\]
where $[i]$ is the equivalence class of $i$ in $\mathbb{Z}_{n+4}/\sim$.
As $\mathbb{Z}_{4m}/\sim$ and $\mathbb{Z}_{4m-4}$ are finite sets, it suffices to show $f(c_q([i]))=c_{m}(f([i]))$ for all $[i]\in \mathbb{Z}_{4m}/\sim$ for $f$ to be a homeomorphism. By \cref{ex:coequalizer}  we have the following calculations:

\begin{itemize}[left=0pt]
\item Suppose $i\in A_1$. Then $c_q([i])=\{[3m-1],[3m],\dots [0],[1],\dots,[m],[m+2]\}$. 
Therefore
\[f(c_q([i]))=\{3m-4,3m-3,\dots ,4m-5,0,1,\dots m-1,m\}=c_{m}(0)=c_{m}(f([i])).\]

\item Suppose $i\in \{2,3,\dots, m-1\}$. Note that \[i-m-1=4m+i-m-1 \textnormal{ mod(4m)}=3m+i-1 \textnormal{ mod(4m)}.\] 
Hence
\[c_q([i])=\{[3m+i-1],[3m+i],\dots ,[0],[1],\dots [i],[i+1],\dots ,[i+m],[i+m+1]\}.\] 
As $i\ge 2$, we have that $3m+i-1\ge 3m+1$.
Therefore
\[f(c_q([i]))=\{3m+i-5,3m+i-4,\dots, i-1,i,\dots, i+m-1 \}=c_{m}(i-1)=c_{m}(f([i])).\]

\item Suppose $i\in A_2$. Then $c_q([i])=\{[4m-1],[0],[1],\dots ,[m],[m+1],\dots ,[2m+2]\}$. 
Therefore 
\[f(c_q([i]))=\{4m-5,0,1,\dots, m-1,m,\dots ,2m-1\}=c_{m}(m-1)=c_{m}(f([i])).\]
\item Suppose $i\in \{m+2,m+3,\dots,2m-1\}$. Then, we have $i-m-1\in \{1,2,\dots, m-2\}$, and $i+m+1\in \{2m+3,2m+4,\dots,3m\}$. Hence $f([i-m-1])=i-m-1-1=i-m-2$. If $i+m+1=3m$, then $f([i+m+1])=f([3m])=3m-3=i+m+1-3=i+m-2$. If $i+m+1\neq 3m$, then $f([i+m+1])=i+m-2$. Since
\[c_q([i])=\{[i-m-1],[i-m],\dots ,[i],[i+1],\dots, [i+m], [i+m+1]\},\] 
we have
\begin{align*}
f(c_q([i]))=&\{i-m-2,i-m-1,\dots ,i-2,i-1,\dots , i+m-3,i+m-2\}=\\
=&c_{m}(i-2)=c_{m}(f([i])).
\end{align*}
\item Suppose $i\in A_3$. Then $c_q([i])=\{[m-1],[m],\dots ,[2m],[2m+1],\dots ,[3m+2]\}$. 
Therefore
\[f(c_q([i]))=\{m-2,m-1,\dots, 2m-2,2m-1, \dots, 3m-2\}=c_{m}(2m-2)=c_{m}(f([i])).\]
\item Suppose $i\in \{2m+2,2m+3,\dots,3m-1\}$. Then, we have $i-m-1\in \{m+1,m+2,\dots ,2m-2\}$, and $i+m+1\in \{3m+3,3m+4,\dots, 4m\}$. Note that 
$4m=0\textnormal{ mod(4m)}$. We have
\[c_q([i])=\{[i-m-1],[i-m],\dots [i],[i+1],\dots ,[i+m],[i+m+1]\}.\] 
If $i-m-1=m+1$, then $f([i-m-1])=m-1=m+1-2=i-m-1-2=i-m-3$. If $i-m-1=2m-2$, then $f([i-m-1])=2m-4=2m-2-2=i-m-1-2=i-m-3$. If $i+m+1=3m+3$, then $f([i+m+1])=3m+3-4=i+m+1-4=i+m-3$. If $i+m+1=0=4m\textnormal{ mod(4m)}$, then $f([i+m+1])=0=4m-4 \textnormal{ mod(4m-4)}= i+m+1-4\textnormal{ mod(4m-4)}=i+m-3\textnormal{ mod(4m-4)}$.
Therefore  
\begin{align*}
f(c_q([i]))=&\{i-m-3,i-m-2,\dots ,i-3,i-2,\dots , i+m-4,i+m-3\}=\\
=&c_{m}(i-3)=c_{m}(f([i])).
\end{align*}
\item Suppose $i\in A_4$. Then 
$c_q([i])=\{[2m-1],[2m],\dots ,[3m],[3m+1],\dots ,[0],[1],[2]\}$.
Therefore 
\[
    f(c_q([i]))=\{2m-3,2m-1,\dots ,3m-3,3m-2,\dots,0,1\}=c_{m}(3m-3)=c_{m}(f([i])).
\]
\item Suppose $i\in \{3m+2,3m+3,\dots ,4m-1\}$. Then, $i-m-1\in \{2m+1,2m+2,\dots ,3m-2\}$, and $i+m+1\in \{3\textnormal{ mod(4m)},4\textnormal{ mod(4m)},\dots ,m\textnormal{ mod(4m)}\}$. If $i-m-1=2m+1$, we have $f([i-m-1])=2m-2=2m+1-3=i-m-1-3=i-m-4$. If $i-m-1=3m-2$, we have $f([i-m-1])=3m-2-3=i-m-1-3=i-m-4$. If $i+m+1=3\textnormal{ mod(4m)}$, that is $i=3m+2$, we have $f([3])=2$ which is a distance $m$ from $f([i])=f([3m+2])=3m-2$ in $\mathbb{Z}_{4m-4}$. If $i+m+1=m\textnormal{ mod(4m)}$, that is $i=4m-1$, we have $f([m])=m-1$ which is a distance $m$ from $f([i])=f([4m-1])=4m-5=i-4$ in $\mathbb{Z}_{4m-4}$. Since 
\[c_q([i])=\{[i-m-1],[i-m],\dots, [i],[i+1],\dots ,[i+m\textnormal{ mod(4m)}],[i+m+1\textnormal{ mod(4m)}]\},\] 
we have
\begin{align*}
    f(c_q([i]))=&\{i-m-4,i-m-3,\dots,i-4,i-3,\dots, 0,1,\dots, i+m-5\textnormal{ mod (4m-4)},\\
    &i+m-4\textnormal{ mod (4m-4)}\}=c_{m}(i-4)=c_{m}(f([i])).
\end{align*}
\end{itemize}
Thus $(\mathbb{Z}_{n+4}/\sim, c_q)$ is homeomorphic to $(\mathbb{Z}_n,c_{m})$ and therefore the result follows.
\end{proof}

\section{Discussion}
\subsection{Summary}
This manuscript is part of a series of works that develop algebraic topology for objects that are not topological spaces but appear often in applications, such as (undirected) graphs. An \emph{undirected graph} is a set with a reflexive and symmetric relation $(X,E)$. Undirected graphs are a full subcategory of $\cat{Cl}$. In particular, the closure of $i\in X$ is the out-neighborhood of $i$ in $(X,E)$, namely $c(i)=\{j\,|\, iEj\}$ and we set $c(A)=\bigcup_{i\in A}c(i)$. Conversely, given a closure space $(X,c)$ with the properties that for all $A\subset X$ we have $c(A)=\bigcup_{a\in A}c(a)$ and $x\in c(y)$ whenever $y\in c(x)$, $(X,c)$ is equivalent to an undirected graph $(X,E)$ by declaring that $xEy$ if and only if $y\in c(x)$.
See \cite{bubenik2024homotopy,dikranjan2013categorical} for more details.

Thus we can think of the roots of unity closure spaces $(\mathbb{Z}_n,c_m)$ as undirected graphs and therefore in this work we have computed singular homology groups of some undirected graphs that are not topological spaces. The definition of fundamental groups and singular homology groups presented here are the natural extension from the familiar definitions for topological spaces. Furthermore, they can be applied to any undirected graph and they provide an alternative perspective compared to many homotopy and homology theories studied in applied topology (see Related works subsection).

\subsection{Difficulties}
The computations in \cref{section:applications} are somewhat partial results; given $n,m\in \mathbb{N}$, and $k\ge 0$, we only derived some of $H_k(\mathbb{Z}_n,c_m)$. For example we have not computed $H_k(\mathbb{Z}_6,c_2)$ for $k\ge 1$. Computing singular homology groups is hard due to the infinite dimensional nature of the chain groups even if the space in question has finitely many points. There are infinitely many continuous maps $|\Delta^{k}|\to (\mathbb{Z}_6,c_2)$ for example. That is why we relied on Mayer-Vietoris long exact sequences and the Hurewicz theorem.

One of the difficulties in applying Mayer-Vietoris type arguments is the requirement to have an interior cover, as emphasized in \cref{example:cover_that_is_not_interior_cover}. If we pick any two distinct points in $(\mathbb{Z}_6,c_2)$, their smallest neighborhood will be $\mathbb{Z}_6$. Therefore, the only way to have a non-trivial interior cover for $(\mathbb{Z}_6,c_2)$ would be to consider all six of the sets $c_2(i)$, $i\in \mathbb{Z}_6$. Therefore we cannot construct any meaningful long exact sequences for a 2-element interior cover in this case.

The other difficulty is requiring the existence of good pairs like we needed in proofs of \cref{theorem:higher_homology_for_n_greater_than_6m,theorem:higher_homology_for_n_greater_than_4m} for example. In particular, $c_2(0)=\{0,1,2,4,5\}$ is the smallest neighborhood of $0$ in $(\mathbb{Z}_6,c_2)$ but we don't know if it is a good pair as the smallest neighborhood of $c_2(0)$ is $\mathbb{Z}_6$ and we don't know if $\mathbb{Z}_6$ deformation retracts to $c_2(0)$. 

Therefore, even though the computations present in this manuscript are partial results, they are about the best we can do with the current tools at our disposal. However, there is a way to compute all the $H_k(\mathbb{Z}_n,c_m)$ (and singular homology groups of any finite closure space in fact), that we outline in the next subsection.

\subsection{Concurrent work}
Singular homology groups of topological spaces (and closure spaces) are difficult to compute directly. Tools like long exact sequences of interior covers and excision are often used to simplify the problem. Another strategy is to find a CW complex or a simplicial complex with the same homotopy type (or weak homotopy type) as the space in question and compute their simplicial or cellular homology. Finite topological spaces can be thought of as preorders and it was shown they have the weak homotopy type of their associated order complex by McCord in \cite{mccord1966singular}.

In a companion manuscript \cite{milicevic2024directed} we show that finite directed graphs (where the relation is no longer required to be symmetric) have the weak homotopy type of their directed clique complexes, thus extending the work in \cite{mccord1966singular} to finite closure spaces. This statement for undirected graphs was shown very recently in \cite{trevino2024graphs} using pseudotopological spaces, a generalization of closure spaces. In classical algebraic topology weak homotopy equivalences induce isomorphisms on homology groups but surprisingly this result was not extended to pseudotopological spaces in \cite{trevino2024graphs}. In \cite{milicevic2024directed} we show that weak homotopy equivalences still induce isomorphisms on singular homology in this generalized setting and thus we can calculate singular homology groups of finite graphs like $(\mathbb{Z}_n,c_m)$ by constructing their clique complexes and calculating their simplicial homology instead. 

\subsection*{Acknowledgments}
The author would like to thank Peter Bubenik for helpful discussions and suggestions, and the anonymous referee for their many helpful recommendations that improved this manuscript.

\printbibliography

\end{document}